\documentclass[11pt]{article}
\author{Damian Sercombe}

\title{The length and depth of real algebraic groups}

\usepackage{amssymb,amsmath,pifont,amsfonts,tkz-berge,amsthm,enumitem}
  \usepackage{times}
\usepackage{mathtools}
  \usepackage[T1]{fontenc}
  \usepackage{textcomp}
 \usepackage{lipsum}
\usepackage{multirow,boldline}
\usepackage[margin=1.2in]{geometry}
\usepackage{makecell}
\usepackage{tikz}\usetikzlibrary{shapes,fit}
\usepackage{scalerel}
\usepackage{relsize}
\usepackage{tikz-cd}
\usetikzlibrary{positioning}
\usepackage{rotating}
\usepackage{dynkin-diagrams}
\usepackage{longtable} 
\usepackage[tableposition=above,width=12in]{caption}
\usepackage{arydshln} 
\usepackage{enumitem}

\let\OLDthebibliography\thebibliography
\renewcommand\thebibliography[1]{
  \OLDthebibliography{#1}
  \setlength{\parskip}{0pt}
  \setlength{\itemsep}{0pt plus 0.3ex}
}

\numberwithin{equation}{section}

\theoremstyle{definition}
\newtheorem{thm}{Theorem}

\newtheorem{Remark}[equation]{Remark}

\newtheorem{Theorem}[equation]{Theorem}
\newtheorem{Proposition}[equation]{Proposition}
\newtheorem{Lemma}[equation]{Lemma}
\newtheorem{Corollary}[equation]{Corollary}

\DeclareMathOperator{\PSL}{PSL}

\DeclareMathOperator{\Sp}{Sp}

\DeclareMathOperator{\Exp}{exp}

\DeclareMathOperator{\Aut}{Aut}

\DeclareMathOperator{\Out}{Out}

\DeclareMathOperator{\Gal}{Gal}
\DeclareMathOperator{\SL}{SL}
\DeclareMathOperator{\U}{U}

\DeclareMathOperator{\SU}{SU}
\DeclareMathOperator{\PSU}{PSU}

\DeclareMathOperator{\diag}{diag}

\DeclareMathOperator{\Cl}{Cl}

\newcommand{\R}{\mathbb{R}}
\newcommand{\Z}{\mathbb{Z}}

\newcommand{\C}{\mathbb{C}}
\newcommand{\T}{\mathbb{T}}

\begin{document}

\pagenumbering{arabic}

\maketitle

\section*{Abstract}
\vspace{-1.5mm}
\noindent Let $G$ be a connected real algebraic group. An unrefinable chain of $G$ is a chain of subgroups $G=G_0>G_1>...>G_t=1$ where each $G_i$ is a maximal connected real subgroup of $G_{i-1}$. The maximal (respectively, minimal) length of such an unrefinable chain is called the length (respectively, depth) of $G$. We give a precise formula for the length of $G$, which generalises results of Burness, Liebeck and Shalev on complex algebraic groups \cite{BLS} and also on compact Lie groups \cite{BLS1}. If $G$ is simple then we bound the depth of $G$ above and below, and in many cases we compute the exact value. In particular, the depth of any simple $G$ is at most $9$.
\vspace{-3mm}
\section{Introduction}
\vspace{-1.5mm}
\noindent Let $G$ be a connected real algebraic group. An \textit{unrefinable} chain of length $t$ of $G$ is a chain of real subgroups $G=G_0>G_1>...>G_t=1$ where each $G_i$ is a maximal connected real subgroup of $G_{i-1}$. The \textit{length} $l(G)$ (resp. \textit{depth} $\lambda(G)$) of $G$ is the maximal (resp. minimal) length of such an unrefinable chain. The corresponding notions for connected complex algebraic groups are denoted by $l_{\C}$ and $\lambda_{\C}$. Let $G(\C)$ denote the complexification of $G$ and let $R(G)$ be the radical of $G$. 

\vspace{2mm}\noindent In this paper we study the length and depth of real algebraic groups. These invariants were first introduced for finite groups in the 1960s (see \cite{BLS3} and the references therein for a comprehensive summary). More recently, length and depth have been introduced and studied for algebraic groups over algebraically closed fields in \cite{BLS} and for compact real Lie groups in \cite{BLS1}. In this paper we generalise the latter results by looking at all real algebraic groups.

\vspace{2mm}\noindent In Theorem \ref{length} we obtain a precise formula for the length of any  connected reductive real algebraic group. In addition, we bound the depth of any simple real algebraic group in Theorem \ref{depth}. To prove Theorems \ref{length} and \ref{depth} we use the exact values for the length and depth of any simple complex algebraic group and any simple compact Lie group as computed in \cite{BLS} and \cite{BLS1} respectively.

\vspace{2mm}\noindent Henceforth let $G$ be a connected reductive real algebraic group. The semisimple quotient $G/R(G)$ is $\R$-isomorphic to the derived subgroup $G'$. We can decompose $G$ as a commuting product $\big(\prod_{i=1}^m G_i\big) \times \mathcal{T}^k$ where each $G_i$ is simple and $\mathcal{T}^k$ is a torus of dimension $k$. The \textit{rank} $r(G)$ of $G$ is the dimension of a maximal torus of $G$ and the \textit{real rank} $r_{\R}(G)$ of $G$ is the dimension of a maximal $\R$-split torus of $G$. The \textit{semisimple rank} and the \textit{semisimple real rank} of $G$ refer to the quantities $r(G')$ and $r_{\R}(G')$ respectively. Up to conjugacy, there exists a unique maximal compact subgroup $K$ of $G$ and a unique compact form $G_c$ of $G(\C)$ (refer to $\S \ref{realforms}$).

\vspace{2mm}\noindent Let $S$ be a maximal $\R$-split torus of $G$ and let $T$ be a maximal real torus of $G$ that contains $S$. Let $W=N_G(T)/T$ be the Weyl group of $G$ with respect to $T$ and let $w_0$ be the longest element of $W$. Let $\Phi$ be the root system of $G$ with respect to $T$, let $\Delta$ be a base of $\Phi$ and let $\Phi^+$ be the corresponding subset of positive roots of $\Phi$. Let $\Delta_0$ be the subset of $\Delta$ that vanishes on $S$. The derived subgroup $C_G(S)'$ is a compact real subgroup of $G$ called the \textit{(semisimple) anisotropic kernel} of $G$. Recall from $\S 35.5$ of \cite{H} that $\Delta_0$ is a base of the root system $\Phi_0 \subset \Phi$ of $C_G(S)'$.  

\vspace{2mm}\noindent As described in $\S 2.3$ of \cite{T}, there is a natural action $*$ of the Galois group $\Gamma = \Gal(\C / \R) \cong \Z_2$ on $\Delta$ that stabilises $\Delta_0$. The orbits of $\Gamma$ in $\Delta \setminus \Delta_0$ are called \textit{distinguished}. The \textit{index} $\mathcal{S}(G)$ of $G$ is the data consisting of $\Delta$, $\Delta_0$ and the $*$-action of $\Gamma$ on $\Delta$. We illustrate $\mathcal{S}(G)$ using a \textit{Tits-Satake diagram}, which is constructed by taking the Dynkin diagram of $G$, blackening each vertex in $\Delta_0$ and linking all of the vertices in each $\Gamma$-orbit of $\Delta$ with a solid gray bar.

\vspace{2mm}\noindent Let $\big\{\mathcal{O}_j \subset \Delta \hspace{0.5mm}\big|\hspace{0.5mm} 1 \leq j \leq r_{\R}(G)\big\}$ be the set of distinguished orbits in $\Delta \setminus \Delta_0$. Associated to any $I:=\Delta \setminus (\cup_j\hspace{0.5mm}\mathcal{O}_j)$ is a parabolic subsystem $\Phi_{I}$ of $\Phi$ with base $I$, a standard real parabolic subgroup $P_I$ of $G$ and a standard real Levi subgroup $L_{I}$ of $G$ (see $\S 21.11$ of \cite{B}). If $I=\Delta_0$ then $\Phi_{\Delta_0}=\Phi_0$, $P_{\Delta_0}$ and $L_{\Delta_0}=C_G(S)$ are called \textit{minimal}.

\vspace{2mm}\noindent Our arguments and results are independent of the choice of isogeny type. Theorem \ref{length} gives a formula for the length of any connected reductive real algebraic group. To obtain the explicit values of this formula we need the results of \cite{BLS1} (see $\S \ref{prelimlength}$) to compute the length of the compact group $C_G(S)'$. 

\begin{thm}\label{length} Let $G$ be a connected reductive real algebraic group. Then $$l(G) = |\Phi^+|-|\Phi_0^+|+r(G)+r_{\R}(G')-r\big(C_G(S)'\big)+l\big(C_G(S)'\big) \hspace{1mm}.$$ In particular, the values of $l(G)$ for $G(\C)$ a simple complex group are given in Tables \ref{classical} and \ref{exceptional}.
\end{thm}

\vspace{0mm}\noindent In Theorem \ref{depth} we bound the depth of any simple real algebraic group $G$. The case where $G$ is compact has already been done in \cite{BLS1}. If $G(\C)$ is a simple complex group then the depth $\lambda_{\C} \big(G(\C)\big)$ has been computed in \cite{BLS}. In particular, $3 \leq \lambda_{\C} \big(G(\C)\big) \leq 6$. The roman numeral notation used for exceptional $G$ is standard in the literature, for example see Figure $6.2$ of \cite{Kn}. 

\begin{thm}\label{depth} Let $G$ be a simple real algebraic group. Then $\lambda_{\C} \big(G(\C)\big)-1 \leq \lambda(G) \leq 9$. Moreover \begin{enumerate}[label=(\roman*)]
\item If $G$ is quasisplit or compact then $\lambda(G)=\lambda_{\C} \big(G(\C)\big)-1$.
\vspace{-1mm}\item For $G$ exceptional $$ \lambda(G) = \begin{cases} 3 & \mbox{if } G=GI,~ FI, ~EV \textnormal{ or } EVIII \\ 4 &\mbox{if } G=FII, ~EI, ~EII,~ EIV \textnormal{ or } EVI \\ 5 & \mbox{if } G=EIII,~ EVII \textnormal{ or } EIX \end{cases}$$ 
\vspace{-3mm}\item For $G$ classical 
\begin{itemize} 
\item if $G=\SL_n(\mathbb{H}) ~(n>1) \textnormal{ or } \mathrm{SO}(2k\hspace{-0.5mm}+\hspace{-0.5mm}1,1) ~(k >3)$ then $\lambda(G)=4$,
\item if $G=\mathrm{SO}^*(2k) ~(k \geq 4)$ then $4 \leq \lambda(G) \leq 6 - \zeta_k$, 
\item if $G=\Sp(p,q) ~(p \geq q>0)$ then $4 \leq \lambda(G) \leq 6-\delta_{pq}-\delta_{1q}$,
\item if $G=\mathrm{SO}(p,q) ~(p \geq q>0)$ then $\lambda(G) \leq 8 - \eta_{pq}$, and 
\item if $G=\SU(p,q) ~(p\geq q>0)$ then $\lambda(G) \leq 9 - \eta_{pq}$
\end{itemize}
where $\zeta_k:=\begin{cases} 1 & \textnormal{if $k \neq 7$ and $k$ is odd } \\ 0 & \textnormal{otherwise }\end{cases}$ and $\eta_{pq}:=\begin{cases} 3 & \textnormal{if $p-q = 3$ or $q=1$ } \\ 2 & \textnormal{if $p-q = 4$ or $q=2$ } \\ 0 & \textnormal{if $q = 7$ and $p\geq 14$ is even } \\ 1 & \textnormal{otherwise }\end{cases}$.
\end{enumerate} The upper bounds of $\lambda(G)$ for $G$ classical are given in Table \ref{classical} and the values of $\lambda(G)$ for $G$ exceptional are given in Table \ref{exceptional}.
\end{thm} 

\vspace{2mm}\noindent I would like to thank my supervisor Prof. Martin Liebeck for introducing me to this field, and for his help and support throughout this research. I also acknowledge the support of an Imperial College PhD Scholarship.

\begin{table}[!htb]\begin{center}
\caption{Length and depth of the classical real algebraic groups}\label{classical}\hspace*{-3cm}\begin{tabular}{| c | c | c | c | c | c | c | c |}    \hline

  $G$ & $r(G)$ & $\!r_{\R}(G)\!$  & $\mathcal{S}(G)$ & $C_G(S)'$ & $l(G)$ & ${K}^{\circ}$ & $\lambda(G)$ \topstrut \\ \hline \hline

\makecell{$\SL_{n+1}(\R)$ \\ $(n \geq 1)$} & $n$ & $n$ & $\dynkin{A}{I}$ & $1$ & $\frac{n}{2}(n\hspace{-0.5mm}+\hspace{-0.5mm}5)$ & $\mathrm{SO}(n\hspace{-0.5mm}+\hspace{-0.5mm}1)$ & \makecell{$2 \textnormal{ if } n=1$ \\ $3 \textnormal{ if } n=2$ \\ $5 \textnormal{ if } n=6$ \\ $\!4 \textnormal{ otherwise}\!$}   \\ \hline
\makecell{$\SL_n(\mathbb{H})$ \\ $(n > 1)$} & $2n\hspace{-0.5mm}-\hspace{-0.5mm}1$ & $\!\!n\hspace{-0.5mm}-\hspace{-0.5mm}1\!\!$ & $\dynkin{A}{II}$ & $\SL(2)^n$ & $2(n^2\hspace{-0.5mm}+\hspace{-0.5mm}n\hspace{-0.5mm}-\hspace{-0.5mm}1)$ & $\Sp(n)$ & $4$ \\ \hline
\makecell{$\SU(p,p)$ \\ $(p > 1)$} & $2p\hspace{-0.5mm}-\hspace{-0.5mm}1$ & $p$ & $\begin{tikzpicture}[baseline=-1.8ex]\dynkin[foldradius=2mm]{A}{IIIb}\end{tikzpicture}$ & $1$ & $2p^2\hspace{-0.5mm}+\hspace{-0.5mm}2p\hspace{-0.5mm}-\hspace{-0.5mm}1$ &$\!\U(p)\!\times\!\SU(p)$ & $4$ \\ \hline
\makecell{$\!\!\SU(p,p\!-\!1)\!\!$ \\ $(p > 1)$} & $2p\hspace{-0.5mm}-\hspace{-0.5mm}2$ & $\!\!p\hspace{-0.5mm}-\hspace{-0.5mm}1\!\!$ & $\!\begin{tikzpicture}[baseline=-1.8ex]\dynkin[foldradius=2mm, ply=2]{A}{oo.oo.oo}\end{tikzpicture}$ & $1$ & $2(p^2\hspace{-0.5mm}-\hspace{-0.5mm}1)$ & $\!\!\U(p)\!\times\! \SU(p\!-\!1)\!\!$ & \makecell{$3 \textnormal{ if } p=2$ \\ $5 \textnormal{ if } p=4$ \\ $\!4 \textnormal{ otherwise}\!$} \\ \hline
\makecell{$\SU(p,q)$ \\ $\!\!(p \hspace{-0.2mm}>\hspace{-0.2mm} q\hspace{-0.2mm}+\hspace{-0.2mm}1)\!\!$} & $\!\!p\hspace{-0.5mm}+\hspace{-0.5mm}q\hspace{-0.5mm}-\hspace{-0.5mm}1\!\!$ & $q$ & $\!\!\!\begin{tikzpicture}[baseline=-1.8ex]\dynkin[foldradius=2mm]{A}{IIIa}\end{tikzpicture}\!\!$ &  $\!\SU(p\!-\!q)\!$ & $2(pq\hspace{-0.5mm}+\hspace{-0.5mm}p\hspace{-0.5mm}-\hspace{-0.5mm}1)$ & $\!\U(p)\!\times\!\SU(q) \!$ & $\!\! \leq  9 \hspace{-0.5mm}-\hspace{-0.5mm} \eta_{pq}\!$ \\ \hline
\makecell{$\mathrm{SO}(p,q)$ \\ $\!\!(p\hspace{-0.2mm}>\hspace{-0.2mm}q\hspace{-0.2mm}+\hspace{-0.2mm}2)\!\!$} & $\left \lfloor{\frac{p+q}{2}}\right \rfloor$ & $q$ & \makecell{$\!\!\dynkin{B}{ooo.**}\!\!$ \\ $\!\!\textnormal{or }\dynkin{D}{oo.***}\!\!$} & $\!\mathrm{SO}(p\!-\!q)\!$ & $\!pq\hspace{-0.6mm}+\hspace{-0.6mm}p\hspace{-0.6mm}-\hspace{-0.6mm}1\hspace{-0.6mm}+\hspace{-0.6mm}\left \lfloor\hspace{-0.5mm}{\frac{p-q}{4}}\hspace{-0.5mm}\right \rfloor\!$ & $\mathrm{SO}(p)\!\times\! \mathrm{SO}(q)$ & $\!\! \leq  8 \hspace{-0.5mm}-\hspace{-0.5mm} \eta_{pq}\!$ \\ \hline
\makecell{$\!\!\mathrm{SO}(p,p\!-\!1)\!\!\hspace{-0.2mm}$ \\ $(p \geq 3)$} & $p\hspace{-0.5mm}-\hspace{-0.5mm}1$ & $\!\!p\hspace{-0.5mm}-\hspace{-0.5mm}1\!\!$ & $\dynkin{B}{ooo.oo}$ & $1$ & $p^2\hspace{-0.5mm}-\hspace{-0.5mm}1$ & $\!\!\mathrm{SO}(p)\!\times\! \mathrm{SO}(p\!-\!1)\!\!$ & \makecell{$3$ if $p\neq4$ \\ $4$ if $p=4$} \\ \hline
\makecell{$\mathrm{SO}(p,p)$ \\ $(p \geq 4)$}  & $p$ & $p$ & $\dynkin{D}{ooo.ooo}$ & $1$ & $p^2\hspace{-0.5mm}+\hspace{-0.5mm}p$ & $\!\mathrm{SO}(p)\!\times\! \mathrm{SO}(p)\!$ & $4$ \\ \hline
\makecell{$\!\!\mathrm{SO}(p,p\!-\!2)\!\!\hspace{-0.2mm}$ \\ $(p \geq 5)$}  & $p\hspace{-0.5mm}-\hspace{-0.5mm}1$ & $\!\!p\hspace{-0.5mm}-\hspace{-0.5mm}2\!\!$ & $\dynkin[foldradius=2mm, ply=2]{D}{ooo.ooo}$ & $1$ & $p^2\hspace{-0.5mm}-\hspace{-0.5mm}p\hspace{-0.5mm}-\hspace{-0.5mm}1$ & $\!\!\mathrm{SO}(p)\!\times\! \mathrm{SO}(p\!-\!2)\!\!$ & $4$  \\ \hline
\makecell{$\mathrm{SO}^*(2k)$ \\ $(k \geq 4)$}  & $k$ & $\left \lfloor{\frac{k}{2}}\right \rfloor$ & \makecell{$\!\dynkin[foldradius=2mm]{D}{IIIb}\!$ \\ $\!\textnormal{or } \dynkin{D}{IIIa}\!$} & $\hspace{-0.2mm}\!\!\big(\hspace{-0.2mm}\!\SU(2)\big)\hspace{-0.5mm}^{\left \lfloor\hspace{-0.4mm}{\frac{k}{2}}\hspace{-0.4mm}\right \rfloor}\!\!\!$ & $k^2\hspace{-0.5mm}+\hspace{-0.5mm}\left \lfloor{\frac{k}{2}}\right \rfloor$ & $\U(k)$ & $\!\! \hspace{-0.5mm}\leq 6 \hspace{-0.5mm}-\hspace{-0.5mm} \zeta_k\!$ \\ \hline
\makecell{$\Sp_{2n}(\R)$ \\ $(n > 1)$}  & $n$ & $n$ & $\dynkin{C}{I}$ & $1$ & $n(n\hspace{-0.5mm}+\hspace{-0.5mm}2)$ & $\U(n)$ & $3$ \\ \hline
\makecell{$\Sp(p,q)$ \\ $(p \geq q)$} & $p\hspace{-0.5mm}+\hspace{-0.5mm}q$ & $q$ & $\!\dynkin{C}{*o*o.**}\!$ & \makecell{$\!\SU(2)^q \times\!$ \\ $\!\Sp(p\!-\!q)\!$} & $\!\!4pq\hspace{-0.6mm}+\hspace{-0.6mm}3p\hspace{-0.6mm}-\hspace{-0.6mm}1\hspace{-0.6mm}+\hspace{-0.6mm}\delta_{pq}\hspace{-0.3mm}\!\!$ & $\Sp(p) \!\times\! \Sp(q)$ & \makecell{$\hspace{-0.5mm}\!\!\leq\! 6\!-\!\delta_{pq} \hspace{-0.4mm}\!-\!\delta_{pq}\!\! \hspace{-0.4mm}$ \\ $\!3 \textnormal{ if } q=0\!$} \\ \hline
  \end{tabular}\hspace*{-3cm}
\end{center}\end{table}

\begin{table}[!htb]\begin{center}
\caption{Length and depth of the exceptional real algebraic groups}\label{exceptional}\begin{tabular}{| c | c | c | c | c | c | c | c |}    \hline

  $G$ & $r(G)$ & $r_{\R}(G)$ & $\mathcal{S}(G)$ & $C_G(S)'$ & $l(G)$ & ${K}^{\circ}$ & $~\lambda(G)$ \topstrut \\ \hline \hline

$GI$ & $2$ & $2$ & $\dynkin{G}{I}$ & $1$ & $10$ & $(A_1)_c^2$ & $3$   \\ \hline
$(G_2)_c$ & $2$ & $0$ & $\dynkin{G}{2}$ & $(G_2)_c$ & $5$ & $(G_2)_c$ & $3$   \\ \hline
$FI$ & $4$ & $4$ & $\dynkin{F}{I}$ & $1$ & $32$ & $(C_3A_1)_c$ & $3$  \\ \hline
$FII$ & $4$ & $1$ & $\dynkin{F}{II}$ & $(B_3)_c$ & $24$ & $(B_4)_c$ & $4$ \\ \hline
$(F_4)_c$ & $4$ & $0$ & $\dynkin{F}{4}$ & $(F_4)_c$ & $11$ & $(F_4)_c$ & $3$  \\ \hline
$EI$ & $6$ & $6$ & $\begin{tikzpicture}[baseline=0.3ex]\dynkin{E}{I}\end{tikzpicture}$ & $1$ & $48$ & $(C_4)_c$ & $4$  \\ \hline
$EII$ & $6$ & $4$ & $\begin{tikzpicture}[baseline=-1.8ex]\dynkin[ply=2,foldradius=2mm]{E}{oooooo}\end{tikzpicture}$ & $1$ & $46$ & $(A_5A_1)_c$ & $4$ \\ \hline
$EIII$ & $6$ & $2$ & $\begin{tikzpicture}[baseline=-1.8ex]\dynkin[ply=2,foldradius=2mm]{E}{oo***o}\end{tikzpicture}$ & $(A_3)_c$ & $41$ & $(D_5)_c\T$ & $ 5$ \\ \hline
$EIV$ & $6$ & $2$ & $\begin{tikzpicture}[baseline=0.3ex]\dynkin{E}{IV}\end{tikzpicture}$ & $(D_4)_c$ & $37$ & $(F_4)_c$ & $4$ \\ \hline
$(E_6)_c$ & $6$ & $0$ & $\begin{tikzpicture}[baseline=-1.8ex]\dynkin[ply=2,foldradius=2mm]{E}{******}\end{tikzpicture}$ & $(E_6)_c$ & $13$ & $(E_6)_c$ & $4$ \\ \hline
$EV$ & $7$ & $7$ & $\begin{tikzpicture}[baseline=0.3ex]\dynkin{E}{V}\end{tikzpicture}$ & $1$ & $77$ & $(A_7)_c$ & $3$ \\ \hline
$EVI$ & $7$ & $4$ & $\begin{tikzpicture}[baseline=0.3ex]\dynkin{E}{VI}\end{tikzpicture}$  & $(A_1)_c^3$ & $74$ & $(D_6A_1)_c$ & $4$ \\ \hline
$EVII$ & $7$ & $3$ & $\begin{tikzpicture}[baseline=0.3ex]\dynkin{E}{VII}\end{tikzpicture}$  & $(D_4)_c$ & $66$ & $(E_6)_c\T$ & $5$ \\ \hline
$(E_7)_c$ & $7$  & $0$ & $\begin{tikzpicture}[baseline=0.3ex]\dynkin{E}{7}\end{tikzpicture}$ & $(E_7)_c$ & $17$ & $(E_7)_c$ & $3$ \\ \hline
$EVIII$ & $8$  & $8$ & $\begin{tikzpicture}[baseline=0.3ex]\dynkin{E}{VIII}\end{tikzpicture}$ & $1$ & $136$ & $(D_8)_c$ & $3$ \\ \hline
$EIX$ & $8$ & $4$ & $\begin{tikzpicture}[baseline=0.3ex]\dynkin{E}{IX}\end{tikzpicture}$ & $(D_4)_c$ & $125$ & $(E_7A_1)_c$ & $5$ \\ \hline
$(E_8)_c$ & $8$  & $0$ & $\begin{tikzpicture}[baseline=0.3ex]\dynkin{E}{8}\end{tikzpicture}$ & $(E_8)_c$ & $20$ & $(E_8)_c$ & $3$ \\ \hline

  \end{tabular}
\end{center}\end{table}

\newpage\section{Preliminaries}\label{preliminaries}

\noindent In $\S \ref{realcomplexalgebraicgroups}$ we introduce and characterise the notion of a real form of a connected reductive complex algebraic group. We also present Komrakov's classification of reductive maximal connected subgroups of real forms of simple complex algebraic groups. In $\S \ref{prelimlength}$ and $\S \ref{prelimdepth}$ respectively we present some results about the length and depth of real and complex algebraic groups.

\subsection{Real and complex algebraic groups}\label{realcomplexalgebraicgroups}

\noindent Let $X$ be a complex algebraic group. A real algebraic group $G$ is a \textit{real form} of $X$ if $G(\C)$ is $\C$-isomorphic to $X$.

\subsubsection{Real forms}\label{realforms}

\noindent In this subsection we let $X$ be a connected reductive complex algebraic group. There is a bijective correspondence (up to conjugacy) between real forms of $X$, holomorphic involutions of $X$ and antiholomorphic involutions of $X$.

\begin{Proposition}[\textit{Problems $2.3.27$, $3.1.9$, $3.1.10$ and Theorem $2.3.6$ of \cite{OV}}]\label{involution} Let $G$ be a real form of $X$. Then there is a unique antiholomorphic involution $\sigma_G$ of $X$ that fixes $G$ pointwise. Conversely, let $\sigma$ be an antiholomorphic involution of $X$. Then the fixed point set $X^{\sigma}$ is a real form of $X$.
\end{Proposition} 

\begin{Theorem}[\textbf{Weyl}, \textit{Theorems $5.2.8$ and $5.2.9$ of \cite{OV}}]\label{Weyl} There exists a real form of $X$ that is compact as a Lie group. Any two compact real forms of $X$ are conjugate.
\end{Theorem}

\noindent Henceforth we fix a compact real form $G_c$ of $X$. Let $\sigma_c$ be the unique antiholomorphic involution of $X$ that satisfies $G_c=X^{\sigma_c}$.

 \begin{Theorem}[\textbf{Cartan}, \textit{Theorems $5.1.4$ and $5.2.3$ of \cite{OV}}]\label{Cartan} Any holomorphic involution of $X$ has a conjugate $\theta$ that commutes with $\sigma_c$. The map $\theta \mapsto \sigma_c \cdot \theta$ defines a bijection from the set of $\Aut(X)$-conjugacy classes of holomorphic involutions of $X$ to the set of $\Aut(X)$-conjugacy classes of antiholomorphic involutions of $X$.
\end{Theorem}

\noindent So let $\theta$ be a holomorphic involution of $X$ that commutes with $\sigma_c$ and denote $H:=X^{\theta}$. Then $\theta$ stabilises the real form $G:=X^{\sigma_c \cdot \theta}$ of $X$ and $G^{\theta}=H^{\sigma_c \cdot \theta}=:K$ is a real form of $H$. We illustrate this in the following commutative diagram. By Theorem $5.3.3$ of \cite{OV}, $K$ is a maximal compact subgroup of $G$ and any maximal compact subgroup of $G$ is conjugate to $K$.

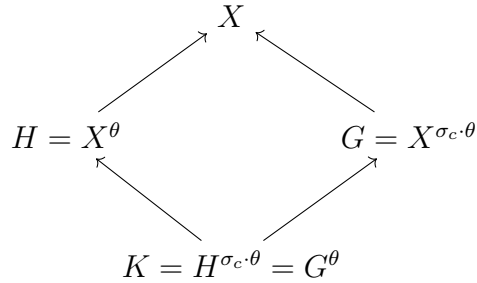
\begin{center}
\begin{tikzpicture}
    \node (E) at (0,0) {\large{$X$}};
    \node[below right=of E] (F) {\large{$G=X^{\sigma_c \cdot \theta }$}};
    \node[below left=of E] (B) {\large{$H=X^{\theta}$}};
    \node (A) at (0,-3.3) {\large{$K=H^{\sigma_c \cdot \theta }=G^{\theta}$}};

    \draw[->] (F)--(E)  ;
    \draw[->] (A)--(F) ;
    \draw[->] (A)--(B) ;
    \draw[->] (B)--(E) ;
\end{tikzpicture}
\end{center}
\noindent If $H$ has maximal rank in $X$ then $G$ is an \textit{inner} form of $X$. Otherwise, $G$ is an \textit{outer} form of $X$.

\begin{Proposition}[\textit{$\S2.2.4$ of \cite{PR}}]\label{splitform} There exists a unique (up to conjugacy) split form $G_s$ of $X$.
\end{Proposition}

\subsubsection{Maximal connected subgroups of simple real algebraic groups}\label{maximalconnsubgroups}

\noindent Any complex algebraic group $X$ can be considered as a real algebraic group $X_{\R}$ of twice the dimension in a process called \textit{realification} (see $\S 2.3.5$ of \cite{OV}). The complexification of $X_{\R}$ is $X^2$. 

\begin{Proposition}[\textit{Theorem $5.1.1$ of \cite{OV}}]\label{simplerealgroups} Let $G$ be a simple real algebraic group. Then either $G=X_{\R}$ for some simple complex group $X$ or $G$ is a real form of a simple complex group. 
\end{Proposition}

\noindent The following result of Komrakov is taken from Tables $3 - 62$ of \cite{K}. However, this source contains a few minor errors, which we correct using Theorem $1$ of \cite{Ta}.

\begin{Theorem}[\textbf{Komrakov}, \cite{K}]\label{Komclass1} Let $G$ be a real form of a simple complex algebraic group such that $G$ is neither split nor compact. Let $M$ be a reductive maximal connected real subgroup of $G$. If $G$ is classical then the possibilities for $M < G$ are listed up to conjugacy in $\Aut(G)$ in Table \ref{classicalcheck} except for the cases where $M(\C)$ is a simple group that acts irreducibly on the natural module of $G(\C)$.  If $G$ is exceptional then the possibilities for $M < G$ are listed up to conjugacy in $\Aut(G)$ in Table \ref{exceptcheck} in $\S \ref{prooflength}$. 

\begin{table}[!htb]\begin{center}
\caption{Reductive maximal connected subgroups of non-split classical real algebraic groups}\label{classicalcheck}\begin{tabular}{| c | c | }    \hline

  $G$ & $M$  \\ \hline \hline

\multirow{3}{*}{\makecell{$\SU(p,q)$ \\ $(p \geq q)$}} & $\SU(p_1,q_1) \hspace{-0.3mm}\times\hspace{-0.3mm} \SU(p_2,q_2) \hspace{-0.3mm}\times\hspace{-0.3mm} \T$ for $p\hspace{-0.2mm}=\hspace{-0.2mm}p_1\hspace{-0.2mm}+\hspace{-0.2mm}p_2$, $q\hspace{-0.2mm}=\hspace{-0.2mm}q_1\hspace{-0.2mm}+\hspace{-0.2mm}q_2$, $(p_1\hspace{-0.2mm}+\hspace{-0.2mm}q_1)(p_2\hspace{-0.2mm}+\hspace{-0.2mm}q_2)\hspace{-0.2mm}\neq \hspace{-0.2mm}0$  \\ 
& $\SU(p_1,q_1) \otimes \SU(p_2,q_2)$ where $p=p_1q_2+p_2q_1$ and $q=p_1p_2+q_1q_2$  \\
& $\SL_n(\C)_{\R}$ where $p=n(n+1)/2$ and $q=n(n-1)/2$ \\ \hline

\multirow{2}{*}{\makecell{$\SL_n(\mathbb{H})$ \\ $(n \geq 2)$}} & $\SL_n(\C)_{\R} \times \T$ \\ 
& $\SL_{n_1}(\R) \otimes \SL_{n_2}(\mathbb{H})$ where $n=n_1n_2$ \\ \hline

\multirow{6}{*}{\makecell{$\mathrm{SO}(p,q)$ \\ $(p > q)$}} & $\mathrm{SO}(p_1,q_1) \times \mathrm{SO}(p_2,q_2)$ for $p=p_1+p_2$, $q=q_1+q_2$, $(p_1+q_1)(p_2+q_2)\neq 0$  \\ 
& $\SU(p/2,q/2) \times \T$ where $p$ and $q$ are even  \\
& $\mathrm{SO}(p_1,q_1) \otimes \mathrm{SO}(p_2,q_2)$ where $p=p_1q_2+p_2q_1$ and $q=p_1p_2+q_1q_2$  \\
& $\Sp(p_1,q_1) \otimes \Sp(p_2,q_2)$ where $p=4(p_1q_2+p_2q_1)$ and $q=4(p_1p_2+q_1q_2)$ \\
& $\mathrm{SO}_n(\C)_{\R}$ where $p=n(n+1)/2$ and $q=n(n-1)/2$  \\ 
& $\Sp_{2n}(\C)_{\R}$ where $p=n(2n+1)$ and $q=n(2n-1)$  \\ \hline

\multirow{5}{*}{\makecell{$\mathrm{SO}^*(2n)$ \\ $(n \geq 3)$}}  & $\mathrm{SO}^*(2n_1) \times \mathrm{SO}^*(2n_2)$ where $n=n_1+n_2$  \\ 
& $\mathrm{SO}_n(\C)_{\R}$  \\
& $\SU(p,q) \times \T$ where $p+q=n$  \\
& $\Sp_{2m}(\R) \otimes \Sp(p,q)$ where $n=2m(p+q)$ \\
& $\mathrm{SO}^*(2m) \otimes \mathrm{SO}(p,q)$ where $n=m(p+q)$ \\ \hline

\multirow{5}{*}{\makecell{$\Sp(p,q)$ \\ $(p \geq q)$}} & $\Sp(p_1,q_1) \times \Sp(p_2,q_2)$ for $p=p_1+p_2$, $q=q_1+q_2$, $(p_1+q_1)(p_2+q_2)\neq 0$  \\ 
& $\SU(p,q) \times \T$  \\ 
& $\Sp(p_1,q_1) \otimes \mathrm{SO}(p_2,q_2)$ where $p=p_1q_2+p_2q_1$ and $q=p_1p_2+q_1q_2$ \\
& $\Sp_{2n_1}(\R) \otimes \mathrm{SO}^*(2n_2)$ for $p=q=n_1n_2$  \\
& $\Sp_{2n}(\C)_{\R}$ where $p=q=n$ \\ \hline
  \end{tabular}
\end{center}\end{table} 
\end{Theorem}

\noindent The following observation is trivial but useful.

\begin{Remark}\label{preserviness} Let $M$ be a proper real subgroup of a real algebraic group $G$. Then $M(\C)$ is a proper complex subgroup of $G(\C)$.
\end{Remark}

\begin{Lemma}\label{maxconn} Let $X$ be a connected reductive complex algebraic group and let $H$ be a connected reductive complex subgroup of $X$. Let $G$ be a real form of $X$ and let $M$ be a real form of $H$ that is contained in $G$. If $H$ is maximal connected in $X$ then $M$ is maximal connected in $G$.
\begin{proof} Let $O$ be a connected proper real subgroup of $G$ that contains $M$. Then ${O(\C)}^{\circ}$ is a connected proper subgroup of $X$ that contains $H$ by Remark \ref{preserviness} and so ${O(\C)}^{\circ}=H$ by maximality. That is, $O$ is a real form of (a finite extension of) $H$ that contains $M$. Hence $O=M$ (again by Remark \ref{preserviness}).
\end{proof} 
\end{Lemma}

\noindent Note that the converse to Lemma \ref{maxconn} does not hold. For example, there is a maximal connected copy of $\PSL_2(\R) \times G_2(\C)_{\R}$ that is contained in the split form of $E_8$ and yet $A_1(G_2)^2<F_4G_2<E_8$ as complex groups.

\begin{Proposition}\label{simplenotsimple} Let $G=X_{\R}$ for some simple complex group $X$. A subgroup $M$ of $G$ is maximal connected if and only if $M$ is maximal parabolic, $M$ is a real form of $X$ or $M=H_{\R}$ for some reductive maximal connected complex subgroup $H$ of $X$. 
\begin{proof} Let $M$ be a non-parabolic maximal connected real subgroup of $G$. Then $M$ is reductive by Corollary $3.3$ of \cite{BT}. It follows from Proposition \ref{involution} that there exists a unique antiholomorphic involution $\sigma$ of $G(\C) \cong X^2$ that satisfies $G=G(\C)^{\sigma}$ and $M=M(\C)^{\sigma}$. 
Let $(x,y) \in G(\C) \cong X^2$. Then $\sigma(x,y)=\big(\sigma_0(y),\sigma_0(x)\big)$ for some antiholomorphic involution $\sigma_0$ of $X$ where $X^{\sigma_0}$ is compact. 
Let $\theta$ be the holomorphic involution of $G(\C)$ that sends $(x,y) \mapsto (y,x)$. Then $G(\C)^{\theta} \cong X$ is diagonally embedded in $G(\C)$. Observe that $\theta$ commutes with $\sigma$ and hence $\theta$ stabilises both $G$ and $M$.

\vspace{2mm}\noindent Let $O$ be a connected reductive $\sigma$-stable proper complex subgroup of $G(\C)$ that contains $M(\C)$. Taking fixed points under $\sigma$ gives us $M \leq O^{\sigma} < G$ and so $O^{\sigma}=M$ by maximality. Then $O=M(\C)$ since $O^{\sigma}$ is a real form of $O$ (this would fail if $O$ was parabolic). That is, $M(\C)$ is maximal among connected reductive $\sigma$-stable complex subgroups of $G(\C)$. Then either $M(\C) \cong X$ is diagonally embedded in $G(\C)$ or $M(\C) \cong H^2$ for some reductive maximal connected complex subgroup $H$ of $X$ where $\sigma$ acts on $M(\C)$ by swapping the two copies of $H$. Hence either $M$ is a real form of $X$ or $M=H_{\R}$.

\vspace{2mm}\noindent Conversely, if $M$ is a real form of $X$ then $M(\C)$ is a diagonally embedded maximal connected subgroup of $G(\C)$ and hence $M$ is maximal connected in $G$ by Lemma \ref{maxconn}. If $M=H_{\R}$ for $H<X$ as stated above then $M$ cannot be contained in a parabolic subgroup of $G$ nor in a real form of $X$. So again $M$ is maximal connected in $G$.
\end{proof}
\end{Proposition}

\begin{Proposition}\label{maxcompact} Let $G$ be a simple real algebraic group and let $K$ be a maximal compact subgroup of $G$. Then ${K}^{\circ}$ is a maximal connected subgroup of $G$.
\begin{proof} Recall the setup from $\S \ref{realforms}$. Let $\sigma_G$ be the unique antiholomorphic involution of $G(\C)$ that satisfies $G=G(\C)^{\sigma_G}$. Let $\theta$ be a holomorphic involution of $G(\C)$ that commutes with $\sigma_G$ and that acts on $G$ with fixed points $G^{\theta}=K$. Denote $H:=G(\C)^{\theta}$.

\vspace{2mm}\noindent We first consider the case where $G=X_{\R}$ for some simple complex group $X$. Recall that $\theta$ acts on $G(\C)=X^2$ by swapping the two copies of $X$  and so $K$ is isomorphic to the compact form of $X$. Hence $K$ is maximal connected in $G$ by Proposition \ref{simplenotsimple}.

\vspace{2mm}\noindent It remains to consider the cases where $G(\C)$ is a simple complex group. Let $\Phi$ be the root system of $G(\C)$ (with respect to some maximal torus), let $\Delta$ be a base of $\Phi$ and let $\alpha_0$ be the longest root of $\Phi$.

 \vspace{2mm}\noindent If ${H}^{\circ}$ is maximal connected in $G$ then the result follows from Lemma \ref{maxconn} since ${K}^{\circ}$ is a real form of ${H}^{\circ}$. Otherwise, we use Table $4.3.1$ of \cite{GLS} (which classifies involutions of simple complex groups) and Theorem $19.1$ of \cite{MT} (which classifies maximal connected subgroups of simple complex groups) to check that ${H}^{\circ}$ is conjugate to a standard Levi subgroup $L_I$ of $G(\C)$ where $I = \Delta \setminus \{\alpha\}$ for some simple root $\alpha$ with a coefficient of $1$ in $\alpha_0$. That is, we can take ${H}^{\circ} < G(\C)$ to be one of the following maximal rank subgroups:
\begin{table}[ht]
\centering 
\begin{tabular}{c | c c c c c c c c c} 
$G(\C)$ & $A_n$ & $B_n$ & $C_n$ & $D_n$ & $D_n$ & $E_6$ & $E_7$ \\ [0.2ex] \hline 
${H}^{\circ}$ & $A_kA_{n-k-1}T^1$ & $B_{n-1}T^1$ & $A_{n-1}T^1$ & $A_{n-1}T^1$ & $D_{n-1}T^1$ & $D_5T^1$ & $E_6T^1$ \topstrut \\ 
\end{tabular}
\end{table}\label{dahdah}

\noindent Assume (for a contradiction) that $M$ is a connected proper real subgroup of $G$ that strictly contains ${K}^{\circ}$. Then ${M(\C)}^{\circ}$ is a connected proper subgroup of maximal rank in $G(\C)$ that strictly contains ${H}^{\circ}$ by Remark \ref{preserviness}. But ${H}^{\circ}$ is maximal among connected reductive subgroups of $G(\C)$ by Corollary $13.7$ and Theorem $13.12$ of \cite{MT} (Borel, de Siebenthal). So $M(\C)$ is conjugate to $P_I$.

\vspace{2mm}\noindent Finally, we observe that ${K}^{\circ}$ contains a compact real torus $T_c$ of dimension $r(G)$ since ${H}^{\circ}$ has maximal rank in $G(\C)$. Then $T_c$ is contained in some Levi subgroup $\mathcal{L} = \mathcal{L}' \times Z(\mathcal{L})^{\circ} \cong M/R_u(M)$ of $M$ where $Z(\mathcal{L})^{\circ} \cong \R^{\times}$ is split. But this is a contradiction as $r(\mathcal{L}')<r(G)$.
\end{proof}
\end{Proposition}

\subsection{Length}\label{prelimlength}

\begin{Theorem}[\textit{Theorem $1$ of \cite{BLS}}]\label{BLSTHM1} Let $G$ be a connected reductive complex algebraic group. Let $B$ be a Borel subgroup of $G$. Then $l_{\C}(G)=\dim(B)+r(G')=|\Phi^+|+r(G)+r(G')$.
\end{Theorem}

\begin{Lemma}[\textit{additivity}]\label{reduction} Let $G$ be a connected real or complex algebraic group with $G/R(G) =\prod_{i=1}^m G_i$ where each $G_i$ is simple. Then $l(G)=\dim\!\big(\hspace{-0.2mm}R(G)\big)+\sum_{i=1}^m l(G_i)$.
\begin{proof} First observe that $l\big(\hspace{-0.2mm}R(G)\big) =\dim\!\big(\hspace{-0.2mm}R(G)\big)$ by Lemma $2.2$ of \cite{BLS} (which holds over both $\R$ and $\C$) 
since $R(G)$ is soluble. The result then follows from Lemma $2.1(ii)$ of \cite{BLS} as $G/R(G)$ is semisimple.
\end{proof} 
\end{Lemma}

\noindent We denote the compact real form of $\Sp_{2n}(\C)$ by $\Sp(n)$ (not $\Sp(2n)$ as in \cite{BLS1}).

\begin{Theorem}[\textit{Theorem $1$ of \cite{BLS1}}]\label{BLSTHM2} The length of each compact simple Lie group $G$ is as follows. 
\begin{table}[ht]
\centering 
\begin{tabular}{c | c c c c c c c c c} 
$G$ & $\SU(n)$ & $\Sp(n)$ & $\mathrm{SO}(n)$ & $G_2$ & $F_4$ & $E_6$ & $E_7$ & $E_8$ \\ [0.2ex] \hline 
$l(G)$ & $2n-2$ & $3n-1$ & $n+\lfloor \frac{n}{4} \rfloor -1$ & $5$ & $11$ & $13$ & $17$ & $20$ \topstrut \\ 
\end{tabular}
\end{table}\label{compactlengths}
\end{Theorem} 

\vspace{-2mm}\noindent In the following lemma we find the lower bound for $l(G)$ as stated in Theorem \ref{length}. We use notation that is taken from the introduction.

\begin{Lemma}\label{lowerbound} Let $G$ be a connected reductive real algebraic group. Then $ l(G) \geq \Lambda_G$ where $$\Lambda_G:= |\Phi^+|-|\Phi_0^+|+r(G)+r_{\R}(G')-r\big(C_G(S)'\big)+l\big(C_G(S)'\big)$$ In particular, if $G$ is split then $C_G(S)'$ is trivial and so $l(G) \geq |\Phi^+|+r(G)+r(G')$. 
\begin{proof} The chain of maximal connected real groups $$G > P_{\Delta \setminus \mathcal{O}_1} > P_{\Delta \setminus (\mathcal{O}_1\cup \mathcal{O}_2)} > ... > P_{\Delta_0}$$ has length $r_{\R}(G')$ and $l(P_{\Delta_0})=\dim\!\big(\hspace{-0.2mm}R(P_{\Delta_0})\big)+l\big(C_G(S)'\big)$ by Lemma \ref{reduction}.

\vspace{2mm}\noindent The $\R$-dimension of any real algebraic group is equal to the $\C$-dimension of its complexification. So $\dim\!\big(C_G(S)'\big)= |\Phi_0| +r\big(C_G(S)'\big)$ by Theorem $8.17(b)$ of \cite{MT} and $\dim(P_{\Delta_0})= |\Phi^+ \cup \Phi_0|+r(G)=|\Phi^+|+ |\Phi^+_0|+r(G)$. Hence $$\dim\!\big(\hspace{-0.2mm}R(P_{\Delta_0})\big)=\dim(P_{\Delta_0})-\dim\!\big(C_G(S)'\big)=|\Phi^+|-|\Phi_0^+|+r(G)-r\big(C_G(S)'\big)$$ and we are done.
\end{proof}
\end{Lemma}

\begin{Proposition}\label{compactsplit} Let $G$ be a non-compact connected reductive real algebraic group where $G'$ is non-trivial. Let $G_c$ (resp. $G_s$) denote the compact (resp. split) form of $G(\C)$. Then $l(G_c) < l(G) \leq l(G_s)=l_{\C}\big(G(\C)\big)$.
\begin{proof} Let $M$ be a maximal connected real subgroup of $G$. Then ${M(\C)}^{\circ}$ is a connected (but not necessarily maximal) proper complex subgroup of $G(\C)$ by Remark \ref{preserviness}. By complexifying an unrefinable chain of $G$ of maximal length we observe that $l(G) \leq l_{\C}\big(G(\C)\big)$. But $l_{\C}\big(G(\C)\big)=|\Phi^+|+r(G)+r(G') \leq l(G_s)$ by Theorem \ref{BLSTHM1} and Lemma \ref{lowerbound} and so $l(G) \leq l(G_s)=l_{\C}\big(G(\C)\big)$. 

\vspace{2mm}\noindent For the lower bound, we write $G$ as a commuting product $\big(\prod_{i=1}^m G_i\big) \times T^k$ where $m \geq 1$, each $G_i$ is simple and $T^k$ is a torus of dimension $k$. Then $l(G) = k+\sum_{i=1}^m l(G_i)$ by Lemma \ref{reduction}. So to show that $ l(G) > l(G_c) $ for any non-compact $G$, it suffices to consider only the cases where $G$ is simple.

\vspace{2mm}\noindent We first consider the cases where $G$ is a non-compact real form of a simple classical complex group. For each case we use Theorem \ref{BLSTHM2}, Lemma \ref{lowerbound} and Table \ref{classical} to check that $ l(G)  \geq \Lambda_G >  l(G_c)$.

\vspace{2mm}\noindent Let $G$ be a non-compact real form of $\SL_n(\C)$. Then $$ l(G) \geq \Lambda_G =\begin{cases} \frac{(n-1)(n+4)}{2} & \mbox{if } G=\SL_n(\R) \\ 2(k^2+k-1) &\mbox{if } G=\SL_{k}(\mathbb{H}) ~(n=2k>2) \\ 2(pq+p-1) +\delta_{pq} & \mbox{if } G=\SU(p,q) ~\Big(\makecell{p \geq q >0, \\  p +q =n}\Big) \end{cases} $$ and $ \Lambda_G > 2n-2 = l\big(\!\SU(n)\big)$. Similarly, let $G$ be a non-compact real form of $\mathrm{SO}_n(\C)$. Then $$ l(G) \geq\Lambda_G =\begin{cases} k^2+\left \lfloor{\frac{k}{2}}\right \rfloor &\mbox{if } G=\mathrm{SO}^*(2k) ~(n=2k \geq 8) \\ pq+p-1+\delta_{pq}+\left \lfloor{\frac{p-q}{4}}\right \rfloor & \mbox{if } G=\mathrm{SO}(p,q) ~\Big(\makecell{p \geq q >0, \\  p+q=n}\Big) \end{cases} $$ and $ \Lambda_G> n+\lfloor \frac{n}{4} \rfloor -1 = l\big(\mathrm{SO}(n)\big)$. If $G$ is a non-compact real form of $\Sp_{2n}(\C)$ then $$ l(G) \geq \Lambda_G =\begin{cases} n(n+2) &\mbox{if } G=\Sp_{2n}(\R) \\ 4pq+3p-1+\delta_{pq} & \mbox{if } G=\Sp(p,q) ~\Big(\makecell{p \geq q >0, \\  p +q =2n}\Big) \end{cases} \hspace{0.3mm}$$ and $ \Lambda_G> 3n-1 =  l\big(\!\Sp(n)\big)$.

\vspace{2mm}\noindent Next, we consider the cases where $G$ is a non-compact real form of an exceptional complex group. Once again we use Theorem \ref{BLSTHM2} and Lemma  \ref{lowerbound} to check that $ l(G) \geq \Lambda_G >  l(G_c)$, where the values of $\Lambda_G$ can be found in Table \ref{exceptional}.

\vspace{2mm}\noindent Finally, let $G=X_{\R}$ where $X$ is a simple complex group with root system $\Phi_X$. Then $G_c=(X_c)^2$ where $X_c$ denotes the compact form of $X$ . It is easy to check that $2l(X_c) \leq l_{\C}(X)+r(X)$ for any simple complex $X$ using Theorems \ref{BLSTHM1} and \ref{BLSTHM2}. Hence $$ l(G) \geq \Lambda_G = 2|\Phi_X^+|+3r(X) > |\Phi_X^+|+3r(X) = l_{\C}(X)+r(X) \geq 2l(X_c) = l(G_c)$$ by Theorem \ref{BLSTHM1} and Lemma \ref{lowerbound}.
\end{proof}
\end{Proposition}

\subsection{Depth}\label{prelimdepth}

\noindent The first lemma is stated in \cite{BLS} over $\C$, but the proof also works over $\R$.

\begin{Lemma}[\textit{Lemma $2.5$ of \cite{BLS}}]\label{lemma2.5} Let $G=NH$ be a connected real or complex algebraic group, where $N$ and $H$ are non-trivial connected subgroups of $G$ and $N \lhd G$. Then $\lambda(G) \geq \lambda(H)+1$.
\end{Lemma}

\begin{Corollary}\label{reductiondepth} Let $G=\big(\prod_{i=1}^m G_i\big) \times T^k$ be a connected reductive real or complex algebraic group where  each $G_i$ is simple and $T^k$ is a torus of dimension $k$. If $m=0$ then $\lambda(G)=k$. If $m \geq 1$ then $$\lambda(G) = \lambda(G') + k \geq \smash{\max\limits_{i=1,...,m}} \big\{\lambda(G_i) \big\}+m-1+k \hspace{1mm}.$$
\begin{proof} It follows immediately from Lemma \ref{lemma2.5} that $\lambda(G) = \lambda(G') + k$. Now assume that $m \geq 1$. We show that $\lambda(G')  \geq \max\limits_{i=1,...,m} \big\{\lambda(G_i) \big\}+m-1$ by induction on $m$. If $m=1$ then we are done. Let $i_0 \leq m$ be a positive integer that satisfies $\lambda(G_{i_0}) \geq \lambda(G_i)$ for all $1 \leq i \leq m$. If $m > 1$ then take some $j \neq i_0$ and so $\lambda(G') \geq \lambda(G'/G_j)+1 \geq \smash{\max\limits_{i=1,...,m}} \big\{\lambda(G_i) \big\}+m-1$ by Lemma \ref{lemma2.5} and the inductive hypothesis.
\end{proof} 
\end{Corollary}

\begin{Corollary}\label{reductiondepthcor} Let $G= (G_0)^m$ where $G_0$ is a simple real or complex algebraic group and $m \geq 1$. Then $\lambda(G)=\lambda(G_0)+m-1$.
\begin{proof} We induct on $m$. If $m=2$ then there exists a diagonally embedded copy of $G_0$ that is maximal connected in $G$ and so $\lambda(G) \leq \lambda(G_0)+1$. If $m>2$ then $\lambda(G) \leq \lambda(G/G_0)+1=\lambda(G_0)+m-1$ by the inductive hypothesis. Then we are done by Corollary \ref{reductiondepth}. 
\end{proof} 
\end{Corollary}

\begin{Theorem}[\textit{Theorem $4$ of \cite{BLS}}]\label{BLSTHM3} Let $G$ be a simple complex algebraic group. Then $$ \lambda_{\C}(G) = \begin{cases} 3 &\mbox{if } G=A_1 \\ 5 & \mbox{if } G=A_r ~(r \geq 3, r \neq 6), B_3, D_r \textnormal{ or } E_6 \\ 6 &\mbox{if } G=A_6 \\ 4 &\mbox{in all other cases} \end{cases}$$
\end{Theorem}

\begin{Theorem}[\textit{Theorem $6$ of \cite{BLS1}}]\label{BLSTHM4} Let $G$ be a compact simple real algebraic group. Then $\lambda(G)=\lambda_{\C} \big(G(\C)\big)-1$.
\end{Theorem}

\section{Proof of Theorem \ref{length}}\label{prooflength}

\noindent Let $G$ be a connected reductive real algebraic group. We already have the lower bound $$l(G) \geq |\Phi^+|-|\Phi_0^+|+r(G)+r_{\R}(G')-r\big(C_G(S)'\big)+l\big(C_G(S)'\big) =: \Lambda_G$$ from Lemma \ref{lowerbound}. So to prove Theorem \ref{length} it suffices to check that $l(M) < \Lambda_G$ for every maximal connected subgroup $M$ of $G$. Our proof is by induction on $l(G)$ and we compute $l(M)$ using Lemma \ref{reduction} (\textit{additivity}), the inductive hypothesis and Tables \ref{classical} and \ref{exceptional}.  

\vspace{2mm}\noindent An outline of the proof is as follows. We first let $M$ be a maximal parabolic subgroup of $G$ and check that $l(M) < \Lambda_G$ by applying the inductive hypothesis to the (reductive) Levi subgroup of $M$. Then we show that it suffices to consider only the cases where $G$ is simple (and neither split nor compact) and $M$ is reductive. Finally, we apply the inductive hypothesis to check that $l(M) < \Lambda_G$ for each of the following cases:
\begin{itemize} 
\vspace{-1.5mm}\item Case (A): $G=X_{\R}$ for some simple complex group $X$.
\vspace{-1.5mm}\item Case (B): $G(\C)=\SL_d(\C)$ for $d \geq 3$ and $M(\C)=\mathrm{SO}_{d}(\C)$ or $\Sp_{d}(\C)$ (if $d$ is even).
\vspace{-1.5mm}\item Case (C): All remaining cases where $G(\C)$ is a simple classical group and $M(\C)$ is a simple group that acts irreducibly on the natural module of $G(\C)$.
\vspace{-1.5mm}\item Case (D): All remaining $M$ in classical $G$ (which are listed in Table \ref{classicalcheck}).
\vspace{-1.5mm}\item Case (E): All remaining $M$ in exceptional $G$ (which are listed in Table \ref{exceptcheck}).
\end{itemize}
\begin{proof} Let $M$ be a maximal connected (real) subgroup of $G$. The complexification ${M(\C)}^{\circ}$ is a connected (but not necessarily maximal) proper complex subgroup of $G(\C)$. If $G$ is split then $l(M)\leq l_{\C}\big(M(\C)\big)<l_{\C}\big(G(\C)\big)=l(G)$ by Proposition \ref{compactsplit}. The case where $G$ is compact has been done in \cite{BLS1}. So we may assume that $G$ is neither split nor compact. In particular, $G$ is not of type $A_1$.

\vspace{2mm}\noindent We first consider the case where $M$ is a parabolic subgroup of $G$. Let $M=P_I$ where $I=\Delta \setminus \mathcal{O}_i$ for some distinguished orbit $\mathcal{O}_i$ of $G$. The anisotropic kernel of $(L_I)'$ is $C_G(S)'$ and the root system is $\Phi_I$. We compute \begin{align*} l(P_I) &=\dim\!\big(\hspace{-0.2mm}R(P_I)\big)+l\big((L_I)'\big)  ~~~~~~~~~~~~~~~~~~~~~~~~~~~~~~~~~~~~~~~~~~~~~~~~~~~~~~~~~~~~~~~~~\hspace{-2.8mm} \textnormal{(Lemma \ref{reduction})} \\ &=|\Phi^+|\hspace{-0.3mm}-\hspace{-0.3mm}|\Phi_I^+|\hspace{-0.3mm}+\hspace{-0.3mm}r(G)\hspace{-0.3mm}-\hspace{-0.3mm}r\big((L_I)'\big)\hspace{-0.3mm}+\hspace{-0.3mm} l\big((L_I)'\big) ~~~~~~~~~~~~~~~~~~~~~~~~~~~~~~~~~\hspace{-0.3mm} \textnormal{(proof of Lemma \ref{lowerbound})} \\
&= |\Phi^+|\hspace{-0.3mm}-\hspace{-0.3mm}|\Phi_0^+|\hspace{-0.3mm}+\hspace{-0.3mm}r(G)\hspace{-0.3mm}+\hspace{-0.3mm}r_{\R}\big((L_I)'\big)\hspace{-0.3mm}-\hspace{-0.3mm} r\big(C_G(S)'\big)\hspace{-0.3mm}+\hspace{-0.3mm}l\big(C_G(S)'\big) ~~~~\hspace{0.6mm} \textnormal{(induction on $l\big((L_I)'\big)$)} \\
&\leq l(G)- r_{\R}(G')+r_{\R}\big((L_I)'\big) ~~~~~~~~~~~~~~~~~~~~~~~~~~~~~~~~~~~~~~~~~~~~~~~~~~~~~~~~~~~\hspace{-3mm} \textnormal{(Lemma \ref{lowerbound})} \\
& \leq l(G)-1 .
 \end{align*} 

\noindent Henceforth we can assume that $M$ is not parabolic. Then $M$ is reductive by Corollary $3.3$ of \cite{BT}. 
Recall that we can decompose $G$ as a commuting product $\big(\prod_{i=1}^m G_i\big) \times T^a$ where each $G_i$ is simple and $T^a$ is a torus of dimension $a$. Since $M$ is maximal connected in $G$, one of the following three possibilities must occur.

\begin{enumerate}[label=(\roman*)]
\vspace{-1mm}\item If $M=\big(\prod_{i=1}^m G_i\big) \times T^{a-1}$ then $l(M)=l(G)-1$ by Lemma \ref{reduction}.
\vspace{-1mm}\item If $M=\big(\prod_{i\neq j,k} G_i\big) \times G^D_{j} \times T^{a} $ for some $j \neq k$ such that $G_j$ is isogenous to $G_k$ (where $G^D_{j} \cong G_j$ is embedded diagonally in $G_j \times G_k$) then $l(M)=l(G)-l(G_j) <l(G)$ by Lemma \ref{reduction}. 
\vspace{-1mm}\item Otherwise, if $M=\big(\prod_{i\neq j} G_i\big) \times M_j  \times T^{a}$ for some $j$ (where $M_j$ is a reductive maximal connected subgroup of $G_j$) then $l(M)=l(G)-l(G_j)+l(M_j)$ again by Lemma \ref{reduction}. 
\end{enumerate} So to show $l(M) < l(G)$ it suffices to consider only the cases where $G$ is a simple real group.

\vspace{3mm}\noindent \textbf{\underline{CASE (A)}:} 

\vspace{2mm}\noindent We first assume that $G=X_{\R}$ for some simple complex group $X$.

\vspace{2mm}\noindent Let $\Phi_X$ be the root system of $X$. Then $\Phi$ is the union of two perpendicular root systems of type $\Phi_X$. Applying Lemma \ref{lowerbound} gives us $l(G) \geq |\Phi^+|+r(G)+r_{\R}(G) = 2|\Phi_X^+|+3r(X)$. Recall that $M$ is a non-parabolic maximal connected subgroup of $G$. By Proposition \ref{simplenotsimple}, either $M=H_{\R}$ for some reductive maximal connected complex subgroup $H$ of $X$ (with root system $\Phi_H$) or $M$ is a real form of $X$. If $M=H_{\R}$ then \begin{equation}\label{yumyum}|\Phi_H^+|+r(H)+r(H')=l_{\C}(H)< l_{\C}(X)=|\Phi_X^+|+2r(X)\end{equation} by Theorem \ref{BLSTHM1} and so \begin{align*} l(M) &=2|\Phi_H^+|+2r(H)+r(H') ~~~~~~~~~~~~~~~~~~\hspace{0.3mm} \textnormal{(inductive hypothesis)} \\ &<2|\Phi_X^+|+3r(X) ~~~~~~~~~~~~~~~~~~~~~~~~~~~~~~~~~ \textnormal{(Equation \ref{yumyum})} \\
&\leq l(G). ~~~~~~~~~~~~~~~~~~~~~~~~~~~~~~~~~~~~~~~~~~~~~~~~~~\hspace{0.1mm} \textnormal{(Lemma \ref{lowerbound})} 
 \end{align*} Similarly, if $M$ is a real form of $X$ then $l(M) \leq l_{\C}(X)=|\Phi_X^+|+2r(X) < 2|\Phi_X^+|+3r(X) \leq l(G)$.

\vspace{2mm}\noindent Henceforth, by Proposition \ref{simplerealgroups}, we can assume that $G$ is a real form of a simple complex group.

\vspace{3mm}\noindent \textbf{\underline{CASE (B)}:} 

\vspace{2mm}\noindent We next consider the cases where $G(\C)=\SL_d(\C)$ for $d \geq 3$ and $M(\C)=\mathrm{SO}_{d}(\C)$ or $\Sp_{d}(\C)$ (if $d$ is even).

\vspace{2mm}\noindent A non-split real form of $\SL_d(\C)$ is isomorphic to either $\SL_k(\mathbb{H})$ (if $d=2k$) or to $\SU(p,q)$ for some $p +q=d$. A real form of $\mathrm{SO}_{d}(\C)$ is isomorphic to either $\mathrm{SO}^*(2k)$ (if $d=2k$) or to $\mathrm{SO}(p,q)$ for some $p +q =d$. Similarly, a real form of $\Sp_{d}(\C)$ is isomorphic to either $\Sp_d(\R)$ or to $\Sp(p,q)$ for some $p +q =d/2$.

\vspace{2mm}\noindent We first assume that $G=\SL_k(\mathbb{H})$ for $d=2k$. Then $l(G) \geq 2(k^2+k-1)$ by Lemma \ref{lowerbound} and Table \ref{classical}. By Proposition \ref{compactsplit}, the length of $M$ is maximised when $M$ is split in which case $l(M)=l_{\C}\big(M(\C)\big)$. Using Theorem \ref{BLSTHM1}, we check that $l_{\C}\big(\mathrm{SO}_{2k}(\C)\big)=k^2+k<2(k^2+k-1)$ and $l_{\C}\big(\!\Sp_{2k}(\C)\big)=k^2+2k<2(k^2+k-1)$. Hence $l(M) \leq l_{\C}\big(M(\C)\big) < l(G)$ for all possible $M <G$, as required.

\vspace{2mm}\noindent Now let $G=\SU(p,q)$ for $p \geq q$ and $p+q \geq 3$. Then $l(G)\geq 2(pq+p-1)+\delta_{pq}$ by Lemma \ref{lowerbound} and Table \ref{classical}. 
In lieu of precise knowledge about which real forms of $M(\C)$ embed in $G$ (which would take some work to establish) we check that $l(M) < l(G)$ for all real forms $M$ of $M(\C)$ that satisfy $r_{\R}(M) \leq r_{\R}(G)=q$.

\vspace{2mm}\noindent If $M=\mathrm{SO}(p',q')$ for $p' \geq q'$ and $p'+q'=p+q$ then $q'=r_{\R}(M) \leq r_{\R}(G)=q$. So \begin{align*} l(M) &=p'q'+p'-1+\delta_{p'q'}+\left \lfloor{(p'-q')/4}\right \rfloor ~~~~~~~~~~~~~~~ \textnormal{(inductive hypothesis, Table \ref{classical})} \\ & < 2(p'q'+p'-1)+\delta_{p'q'} \\ &\leq 2(pq+p-1)+\delta_{pq} \\ & \leq l(G). ~~~~~~~~~~~~~~~~~~~~~~~~~~~~~~~~~~~~~~~~~~~~\hspace{-0.5mm}~~~~~~~~~~~~~~~~~~~~~~~  \textnormal{(Lemma \ref{lowerbound}, Table \ref{classical})} \end{align*} 

\vspace{1mm}\noindent Next let $M=\Sp(p',q')$ for $p' \geq q'$ and $2(p'+q')=p+q$. Then $q' = r_{\R}(M) \leq q$, $p \leq 2p'+q' $ and $$l(M)=4p'q'+3p'-1+\delta_{p'q'}< 2\big((2p'+q')q'+(2p'+q')-1\big)+\delta_{p'q'} \leq 2(pq+p-1)+\delta_{pq} \leq l(G)$$ by the inductive hypothesis and Lemma \ref{lowerbound}.

\vspace{3mm}\noindent If we take $M=\Sp_{2k}(\R)$ for $p+q=2k$ then $k = r_{\R}(M) \leq q$ and so $p=q=k$. Then $$l(M)=k^2+2k < 2k^2+2k-1 \leq l(G)$$ by the inductive hypothesis and Lemma \ref{lowerbound}.

\vspace{3mm}\noindent Finally, if $M=\mathrm{SO}^*(2k)$ for $p+q=2k$ then $\left \lfloor{k/2}\right \rfloor = r_{\R}(M)\leq q$, $p \leq \left \lceil{3k/2}\right \rceil$ and $$l(M)=k^2+\left \lfloor{k/2}\right \rfloor < 2\big(\!\left \lceil{3k/2}\right \rceil \left \lfloor{k/2}\right \rfloor +\left \lceil{3k/2}\right \rceil-1 \big) \leq 2(pq+p-1)+\delta_{pq} \leq l(G)$$ by the inductive hypothesis and Lemma \ref{lowerbound}. 

\vspace{4mm}\noindent \textbf{\underline{CASE (C)}:} 

\vspace{2mm}\noindent Now consider the (remaining) cases where $G(\C)$ is a classical group and $M(\C)$ is a simple group that acts irreducibly on the natural module of $G(\C)$.

\vspace{2mm}\noindent Let $V$ be a complex vector space of dimension $d>1$ equipped with either the zero form or a non-degenerate bilinear form (symmetric or skew-symmetric). Let $G(\C)=\Cl(V)$ be the group of isometries of $V$ with determinant $1$. Let $M(\C)$ be a connected simple complex group that acts irreducibly on $V$ but that is not isomorphic to either $\mathrm{SO}_d(\C)$, $\Sp_d(\C)$ or to $\SL_d(\C)$. Let $\lambda_V$ be the highest weight of $V$ as an irreducible $M(\C)$-module. For any given $M(\C)$, Lemmas $3.5$ and $3.6$ of \cite{BLS1} give a lower bound $N =N\big(M(\C)\big)$ for the dimension $d$ of $V$.

\begin{table}[!htb]\begin{center}\caption{The values of $N=N\big(M(\C)\big)$ for $M(\C)$ a simple complex algebraic group}\label{replengths}\begin{tabular}{c || c c c c | c c | c c | c | c | c | c | c } 
\multicolumn{1}{c||}{} & \multicolumn{4}{c|}{$\SL_k(\C)$} & \multicolumn{2}{c|}{$\Sp_{2k}(\C)$} & \multicolumn{2}{c|}{$\!\mathrm{SO}_k(\C)\!$} & \multicolumn{1}{c|}{} & \multicolumn{1}{c|}{} & \multicolumn{1}{c|}{} & \multicolumn{1}{c|}{} & \multicolumn{1}{c}{} \\
\makecell{$\hspace{-0.3mm}\!M(\C)\!\hspace{-0.3mm}$ \\ ~} & \makecell{~ \\ $\!2\!$} & \makecell{$\!\!\!k \hspace{-0.5mm} =\!\!\!$\\ $\!3\!$} & \makecell{~ \\ $\!4\!$} & $\!k>4\!$ & $\! \hspace{-0.3mm}k\hspace{-0.5mm}=\hspace{-0.5mm}2\!\hspace{-0.3mm}$ & $\!k>2\!$ & \makecell{$\!7 \hspace{-0.5mm}\leq\hspace{-0.5mm} k \hspace{-0.5mm}\leq\hspace{-0.5mm} 14\!$ \\ $\!k \neq 8\!$}  & \makecell{$\!k=8,\!$ \\ $\!k>14\!$} & \makecell{$\!G_2\!$ \\ ~} & \makecell{$\!F_4\!$ \\ ~} & \makecell{$\!E_6\!$ \\ ~} & \makecell{$\!E_7\!$ \\ ~} & \makecell{$\!E_8\!$ \\ ~} \\ [0.2ex] \hline 
$\!N\!$ & $\!4\!$ & $\!6\!$ & $\!10\!$ & $\!\frac{k}{2}(k\hspace{-0.5mm}-\hspace{-0.5mm}1)\!$ & $\!10\!$ & $\!\frac{k}{2}(k\hspace{-0.5mm}-\hspace{-0.5mm}1)\hspace{-0.5mm}-\hspace{-0.5mm}1\!$ & $2^{\left \lfloor{\frac{k-1}{2}}\right \rfloor}$ & $\!\frac{k}{2}(k\hspace{-0.5mm}-\hspace{-0.5mm}1)\!$ & $\!7\!$ & $\!26\!$ & $\!27\!$ & $\!56\!$ & $\!248\!$ \topstrut \\ 
\end{tabular}\end{center}\end{table}
\noindent Assume (for a contradiction) that $l(M) \geq l(G)$. For a given $M(\C)$, by Proposition \ref{compactsplit}, the length of $M$ is maximised when $M$ is split in which case $l(M)=l_{\C}\big(M(\C)\big)$. Similarly, by Proposition \ref{compactsplit} and Theorem \ref{BLSTHM2}, the length of $G$ is minimised when $G=\mathrm{SO}(N)$ in which case $l(G)=N+\left \lfloor{N/4}\right \rfloor-1$. So we can exclude all $M(\C)$ that satisfy the inequality $l_{\C}\big(M(\C)\big) <  N+\left \lfloor{N/4}\right \rfloor-1$. For example, if $M(\C)=E_8(\C)$ then $N(M)=248$ and $l(M) \leq l_{\C}\big(E_8(\C)\big)=136 < 248 + \left \lfloor{248/4}\right \rfloor-1=309 \leq l(G)$ using Theorem \ref{BLSTHM1}. Using Table \ref{replengths}, this argument excludes the cases where $M(\C)$ is isogenous to $\SL_n(\C)$ for $n \geq 17$, $\Sp_{n}(\C)$ for $n \geq 4$, $\mathrm{SO}_n(\C)$ for $n\neq 7,9,10,12$ or $E_8(\C)$.

\vspace{2mm}\noindent We first assume that $G(\C)=\mathrm{SO}_d(\C)$ or $\Sp_d(\C)$ 
which, by Lemma $78$ of \cite{S}, occurs if and only if $\lambda_V=-w_0 \lambda_V$. 
By Tables $6.6 - 6.53$ of \cite{Lu}, if $\lambda_V=-w_0 \lambda_V$ and $l_{\C}\big(M(\C)\big) \geq  d+\left \lfloor{d/4}\right \rfloor-1 = l\big(\mathrm{SO}(d)\big)$ then $(M(\C),d)$ must be one of $(A_5,20)$, $(B_3,8)$, $(B_4,16)$, $(D_6,32)$, $(G_2,7)$, $(F_4,26)$ or $(E_7,56)$. If $(M(\C),d)$ is $(A_5,20)$, $(D_6,32)$ or $(E_7,56)$ then $G(\C)=\Sp_d(\C)$ 
 by Lemma $79$ of \cite{S}. We can exclude each of these three possibilities using Theorems \ref{BLSTHM1}, \ref{BLSTHM2} and Proposition \ref{compactsplit} since $l_{\C}\big(\hspace{-0.4mm}A_5(\C)\big) = 25 < 29=l\big(\!\Sp(10)\big)$, $l_{\C}\big(D_6(\C)\big) = 42 < 47=l\big(\!\Sp(16)\big)$ and $l_{\C}\big(E_7(\C)\big) = 77 < 83=l\big(\!\Sp(28)\big)$. If $(M(\C),d)$ is $(B_3,8)$, $(B_4,16)$, $(G_2,7)$ or $(F_4,26)$ then $G(\C)=\mathrm{SO}_d(\C)$ 
by Lemma $79$ of \cite{S}. If $G$ is not compact then $G$ is isomorphic to either $\mathrm{SO}^*(d)$ (if $d$ is even) or to $\mathrm{SO}(p,q)$ for some $p \geq q >0$ satisfying $p +q =d$. Observe that $$ l(G) \geq \min \hspace{-0.5mm}\big\{ (d/2)^2+\left \lfloor{d/4}\right \rfloor, pq+p-1 +\left \lfloor{(p-q)/4}\right \rfloor \big\} \geq 2d-3$$ by Lemma \ref{lowerbound}. If $(M(\C),d)$ is $(B_4,16)$, $(G_2,7)$ or $(F_4,26)$ then we use Theorem \ref{BLSTHM1} to check that $l_{\C}\big(M(\C)\big) < 2d-3$. Hence $l(M) \leq l_{\C}\big(M(\C)\big) < l(G)$ for all possible $M < G$ by Proposition \ref{compactsplit}, as required. It remains to consider the case where $(M(\C),d)=(B_3,8)$. If $M=\mathrm{SO}(6,1)$ or $\mathrm{SO}(7)$ then we check that $l(M) < 2d-3=13$ using the inductive hypothesis and Table \ref{classical}. If $M=\mathrm{SO}(4,3)$ or $\mathrm{SO}(5,2)$ then $2 \leq r_{\R}(M) \leq r_{\R}(G)$ and so $G$ must be either $\mathrm{SO}(6,2)$, $\mathrm{SO}(5,3)$ or $\mathrm{SO}(4,4)$, but then $l(M) \leq  l_{\C}\big(M(\C)\big) =15 < l(G)$ by Theorem \ref{BLSTHM1} and Lemma \ref{lowerbound}.

\vspace{2mm}\noindent We now assume that $G(\C)=\SL_d(\C)$ 
and that $M(\C)$ satisfies the tightened inequality $l_{\C}\big(M(\C)\big) \geq 2d-2 = l\big(\!\SU(d)\big) $. By Tables $6.6 - 6.53$ of \cite{Lu}, the only possibility for $(M(\C),d)$ is $(A_4,10)$. We use Theorems \ref{BLSTHM1}, \ref{BLSTHM2} and Proposition \ref{compactsplit} (since $G$ is not compact) to check that $l_{\C}\big(\hspace{-0.4mm}A_4(\C)\big)=18=l\big(\!\SU(10)\big)< l(G)$. We have our contradiction.

\vspace{8mm}\noindent \textbf{\underline{CASE (D)}:} 

\vspace{2mm}\noindent In $(i)$ to $(v)$ we consider the remaining non-parabolic maximal connected real subgroups $M$ of classical real $G$ (which are listed in Table \ref{classicalcheck}) and check that $l(M) < l(G)$. Recall that Lemma \ref{lowerbound} gives a lower bound for $l(G)$ and that we can compute $l(M)$ using Lemma \ref{reduction}, the inductive hypothesis and Tables \ref{classical} and \ref{exceptional}.

\vspace{2mm}\noindent \underline{$(i)$}: $G=\SU(p,q)$ ($p \geq q$, $p+q \geq 2$), where $l(G) \geq 2(pq\hspace{-0.2mm}+\hspace{-0.2mm}p\hspace{-0.2mm}-\hspace{-0.2mm}1) \hspace{-0.2mm}+\hspace{-0.2mm}\delta_{pq}$ by Lemma \ref{lowerbound} and Table \ref{classical}.


\vspace{2mm}\noindent Let $M=\SU(p_1,q_1) \times \SU(p_2,q_2) \times \T$ where $p=p_1+p_2$, $q=q_1+q_2$ and $(p_1+q_1)(p_2+q_2)\neq 0$. We compute \begin{align*} l(M) &=l\big(\!\SU(p_1,q_1)\big) +l\big(\!\SU(p_2,q_2) \big)+l(\T) ~~~~~~~~~~~~~~~~~~~~~~~~~~~~~~~ \textnormal{(Lemma \ref{reduction}, \textit{additivity})}\\ &=2(p_1q_1\hspace{-0.2mm}+p_1\hspace{-0.2mm}-1)\hspace{-0.2mm}+\hspace{-0.2mm}\delta_{p_1q_1}\hspace{-0.6mm}+\hspace{-0.2mm}2(p_2q_2\hspace{-0.2mm}+\hspace{-0.2mm}p_2\hspace{-0.2mm}-\hspace{-0.2mm}1)\hspace{-0.2mm}+\hspace{-0.2mm}\delta_{p_2q_2}\hspace{-0.6mm}+\hspace{-0.2mm}1 \hspace{-0.2mm}~~~ \textnormal{(inductive hypothesis, Table \ref{classical})} \\ & < 2(pq+p-1) +\delta_{pq}  \leq l(G) \hspace{-0.2mm}~~~~~~~~~~~~~~~~~~~~~~~~~~~~~~~~~~~~~~~~~~~~~~~ \textnormal{(Lemma \ref{lowerbound}, Table \ref{classical})} \end{align*}

\noindent Similarly, let $M=\SU(p_1,q_1) \otimes \SU(p_2,q_2)$ where $p=p_1q_2+p_2q_1$ and $q=p_1p_2+q_1q_2$ for non-negative integers $p_1 \geq q_1$ and $q_2 \geq p_2$ that satisfy $p_1+q_1 \geq 2$ and $p_2+q_2 \geq 2$. Then $l(M)=2(p_1q_1+p_1-1)+\delta_{p_1q_1}+2(p_2q_2+q_2-1)+\delta_{p_2q_2} < 2(pq+p-1) +\delta_{pq}$.

\vspace{2mm}\noindent Finally, let $M=\SL_n(\C)_{\R}$ where $p=n(n+1)/2$, $q=n(n-1)/2$ and $n>1$. Then $l(M)=(n+3)(n-1)<
2(pq+p-1)$.

\vspace{2mm}\noindent \underline{$(ii)$}: $G=\SL_n(\mathbb{H})$ ($n \geq 2$), where $l(G) \geq 2(n^2+n-1)$ by Lemma \ref{lowerbound} and Table \ref{classical}.

\vspace{2mm}\noindent If $M=\SL_n(\C)_{\R} \times \T$ then $l(M)=(n+3)(n-1)+1< 2(n^2+n-1)$.

\vspace{2mm}\noindent Let $M=\SL_{n_1}(\R) \otimes \SL_{n_2}(\mathbb{H})$ where $n=n_1n_2$ and $n_1>1$. Then $l(M)=(n_1-1)(n_1+4)/2+2(n_2^2+n_2-1)<  2(n^2+n-1)$.

\vspace{2mm}\noindent \underline{$(iii)$}: $G=\mathrm{SO}(p,q)$ ($p \hspace{-0.2mm}>\hspace{-0.2mm} q$, $p\hspace{-0.3mm}+\hspace{-0.3mm}q\hspace{-0.2mm}\geq\hspace{-0.2mm} 5$), where $l(G) \hspace{-0.2mm}\geq\hspace{-0.2mm} pq\hspace{-0.3mm}+\hspace{-0.3mm}p\hspace{-0.3mm}-\hspace{-0.3mm}1\hspace{-0.3mm}+\hspace{-0.3mm}\left \lfloor\hspace{-0.3mm}{\frac{p-q}{4}}\hspace{-0.3mm}\right \rfloor$ by Lemma \ref{lowerbound} and Table \ref{classical}.

\vspace{2mm}\noindent First we let $M=\mathrm{SO}(p_1,q_1) \times \mathrm{SO}(p_2,q_2)$ where $p=p_1+p_2$, $q=q_1+q_2$, $p_1> q_1$ and $p_2+q_2\neq 0$. We compute $l(M)=p_1q_1+p_1-1+\left \lfloor{(p_1-q_1)/4}\right \rfloor+p_2q_2+\max\{p_2,q_2\}-1+\delta_{p_2q_2}+\left \lfloor{|p_2-q_2|/4}\right \rfloor < pq+p-1 +\left \lfloor{(p-q)/4}\right \rfloor$ (this holds even if $p_2 < q_2$).

\vspace{2mm}\noindent If $M=\SU(p/2,q/2) \times \T$ where $p$ and $q$ are both even and $p+q \geq 4$ then $l(M)=(pq)/2+p-2+1 < pq+p-1 +\left \lfloor{(p-q)/4}\right \rfloor$.

\vspace{2mm}\noindent Now consider $M=\mathrm{SO}(p_1,q_1) \otimes \mathrm{SO}(p_2,q_2)$ where $p=p_1q_2+p_2q_1$ and $q=p_1p_2+q_1q_2$ for non-negative integers $p_1 > q_1$ and $q_2 > p_2$ satisfying $p_1 \geq 2$ and $q_2 \geq 2$. Then $l(M)=p_1q_1+p_1-1+\left \lfloor{(p_1-q_1)/4}\right \rfloor+p_2q_2+q_2-1+\left \lfloor{(q_2-p_2)/4}\right \rfloor < pq+p-1 +\left \lfloor{(p_1-q_1)(q_2-p_2)/4}\right \rfloor = pq+p-1+\left \lfloor{(p-q)/4}\right \rfloor$.

\vspace{2mm}\noindent Next, let $M=\Sp(p_1,q_1) \otimes \Sp(p_2,q_2)$ where $p=4(p_1q_2+p_2q_1)$ and $q=4(p_1p_2+q_1q_2)$ for non-negative integers $p_1 > q_1$ and $q_2 > p_2$. Then $l(M)=(4p_1q_1+3p_1-1)+(4p_2q_2+3q_2-1) \leq 4\big(p_1q_1(p_2^2+q_2^2)+p_2q_2(p_1^2+q_1^2)\big)+(4p_1q_2+2)-2 < pq+p+\left \lfloor{(p-q)/4}\right \rfloor-1$.

\vspace{2mm}\noindent Let $M=\mathrm{SO}_n(\C)_{\R}$ where $p=n(n+1)/2$, $q=n(n-1)/2$ and $n>1$. Then $l(M) \leq n(n+1)/2<(n^4+n^2+2n-4)/4+\left \lfloor{n/4}\right \rfloor=pq+p-1+\left \lfloor{(p-q)/4}\right \rfloor$.

\vspace{2mm}\noindent Finally, let $M=\Sp_{2n}(\C)_{\R}$ where $p=n(2n+1)$ and $q=n(2n-1)$. Then $l(M)=n(2n+3)<4n^4+3n^2+n-1+\left \lfloor{n/2}\right \rfloor=pq+p-1+\left \lfloor{(p-q)/4}\right \rfloor$.

\vspace{2mm}\noindent \underline{$(iv)$}: $G=\mathrm{SO}^*(2n)$ ($n \geq 3$), where $l(G) \geq n^2+\left \lfloor{n/2}\right \rfloor$ by Lemma \ref{lowerbound} and Table \ref{classical}.

\vspace{2mm}\noindent First we take $M=\mathrm{SO}^*(2n_1) \times \mathrm{SO}^*(2n_2)$ where $n=n_1+n_2$ for positive integers $n_1$ and $n_2$. We compute $l(M)=n_1^2+\left \lfloor{n_1/2}\right \rfloor + n_2^2+\left \lfloor{n_2/2}\right \rfloor< n^2+\left \lfloor{n/2}\right \rfloor$.

\vspace{2mm}\noindent If $M=\mathrm{SO}_n(\C)_{\R}$ for $n>1$ then $l(M) \leq n(n+1)/2<n^2+\left \lfloor{n/2}\right \rfloor$.

\vspace{2mm}\noindent Now if $M=\SU(p,q) \times \T$ where $p+q=n$ then $l(M) \leq 2(pq+p) < n^2+\left \lfloor{n/2}\right \rfloor$.

\vspace{2mm}\noindent Let $M=\Sp_{2m}(\R) \otimes \Sp(p,q)$ where $n=2m(p+q)$ and $p \geq q$. Then $l(M)=m(m+2)+4pq+3p-1+\delta_{pq} \leq (2m^2+1) +4pq+2p^2+\delta_{pq} \leq 2\big(m^2+(p+q)^2\big)+1+\delta_{pq}< n^2+\left \lfloor{n/2}\right \rfloor$.

\vspace{2mm}\noindent Finally, let $M=\mathrm{SO}^*(2m) \otimes \mathrm{SO}(p,q)$ where $n=m(p+q)$ and $p+q \geq 2$. Then $l(M)=m^2+\left \lfloor{m/2}\right \rfloor + pq+p-1+\delta_{pq}+\left \lfloor{(p-q)/4}\right \rfloor < (m^2+(p+q)^2-1) +\left \lfloor{n/2}\right \rfloor \leq n^2+\left \lfloor{n/2}\right \rfloor$.

\vspace{2mm}\noindent \underline{$(v)$}: $G=\Sp(p,q)$ ($p \geq q$), where $l(G) \geq 4pq+3p-1+\delta_{pq}$ by Lemma \ref{lowerbound} and Table \ref{classical}.

\vspace{2mm}\noindent First let $M=\Sp(p_1,q_1) \times \Sp(p_2,q_2)$ where $p=p_1+p_2$, $q=q_1+q_2$ and $(p_1+q_1)(p_2+q_2)\neq 0$. We compute $l(M)=(4p_1q_1+3p_1-1+\delta_{p_1q_1})+(4p_2q_2+3p_2-1+\delta_{p_2q_2})< 4pq+3p-1+\delta_{pq}$.

\vspace{2mm}\noindent If $M=\SU(p,q) \times \T$ then $l(M) = 2(pq+p-1)+\delta_{pq}+1 < 4pq+3p-1+\delta_{pq}$.

\vspace{2mm}\noindent Now let $M=\Sp(p_1,q_1) \otimes \mathrm{SO}(p_2,q_2)$ where $p=p_1q_2+p_2q_1$ and $q=p_1p_2+q_1q_2$ for non-negative integers $p_1 \geq q_1$ and $q_2 \geq p_2$ satisfying $q_2 \geq 2$. Then $l(M)=4p_1q_1+3p_1-1+\delta_{p_1q_1}+p_2q_2+q_2-1+\delta_{p_2q_2}+\left \lfloor{(q_2-p_2)/4}\right \rfloor < 4p_1q_1+p_2q_2+3p_1-2+\delta_{p_1q_1}+\delta_{p_2q_2}+\big(3q_2-2)\leq 4(p_1q_1+p_2q_2)+3(p_1+q_2-1)-1+\delta_{p_1q_1} \leq 4pq+3p-1+\delta_{pq}$.

\vspace{2mm}\noindent Next, let $M=\Sp_{2n_1}(\R) \otimes \mathrm{SO}^*(2n_2)$ for $p=q=n_1n_2$. Then $l(M)=n_1(n_1+2)+n_2^2+\left \lfloor{n_2/2}\right \rfloor < 4n_1^2n_2^2+3n_1n_2=4p^2+3p$.

\vspace{2mm}\noindent Finally, let $M=\Sp_{2n}(\C)_{\R}$ where $p=q=n$. Then $l(M)=n(2n+3)<4p^2+3p$.

\vspace{4mm}\noindent \textbf{\underline{CASE (E)}:}

\vspace{2mm}\noindent In Table \ref{exceptcheck} we consider the remaining non-parabolic maximal connected real subgroups $M$ of exceptional real $G$ (which are taken from Tables $4-62$ of \cite{K}) and check that $l(M) < l(G)$. Recall that Lemma \ref{lowerbound} gives a lower bound for $l(G)$ and that we can compute $l(M)$ using Lemma \ref{reduction}, the inductive hypothesis and Tables \ref{classical} and \ref{exceptional}. 

\vspace{2mm}\noindent This completes the proof of Theorem \ref{length}. 
\end{proof}

\newgeometry{top=1.2cm,bottom=1.2cm,left=0cm, right=0cm,foot=5mm}
\begin{longtable}{| c | c | c | c | c | c |}\caption{Reductive maximal connected subgroups of non-split, non-compact exceptional real groups}\label{exceptcheck} \\ \hline

  $G$ & $\mathcal{S}(G)$ & $l(G) \geq$ & $\mathcal{S}(M)$ & $l(M)$ \\ \hline \hline
\endfirsthead \multicolumn{6}{|c|}%
{{\bfseries \tablename\ \thetable{} -- continued from previous page}} \\ \hline
  $G$ & $\mathcal{S}(G)$ & $\!l(G) \!\geq\!$ & $\mathcal{S}(M)$ & $l(M)$ \\ \hline \hline \endhead
\hline \multicolumn{6}{|c|}{{Continued on next page}} \\ \hline \endfoot
 \endlastfoot

\multirow{5}{*}{$FII$} & \multirow{5}{*}{$\dynkin{F}{II}$} & \multirow{5}{*}{$24$} & $\dynkin{B}{o***}$ & $16$ \\ 
&&& $\dynkin{C}{*o*}\times\dynkin{A}{*}$ & $15$ \\
&&& $\dynkin{B}{****}$ & $10$ \\
&&& $\dynkin[ply=2,rotate=90,foldradius=2mm]{A}{oo} \times \dynkin[ply=2,rotate=90,foldradius=2mm]{A}{**}$ & $10$\\
&&& $\dynkin{G}{2} \times \dynkin{A}{o}$ & $8$ \\ \hline

\multirow{6.3}{*}{$EIV$} & \multirow{6.3}{*}{$\dynkin{E}{IV}$} & \multirow{6.3}{*}{$37$} & $\dynkin{A}{*o*o*}\times \dynkin{A}{*}$ & $24$ \\
&&& $\dynkin{F}{II}$ & $24$ \\
&&& $\dynkin{C}{*o**}$ & $20$ \\
&&& $\begin{tikzpicture}[baseline=0.35ex] 
\dynkin[name=upper]{A}{oo}
\node (current) at ($(upper root 1)+(-0,3.7mm)$) {};
\dynkin[at=(current),name=lower]{A}{oo}
\begin{scope}[on background layer]
\foreach \i in {1,2}%
{%
\draw[/Dynkin diagram/foldStyle]
($(upper root \i)$) -- ($(lower
root \i)$);%
}%
\end{scope}
\end{tikzpicture} \times \dynkin[ply=2,foldradius = 2mm,rotate=90]{A}{**}$ & $16$ \\ 
&&& $\dynkin{G}{2} \times \dynkin{A}{oo}$ & $12$ \\ 
&&& $\dynkin{F}{4}$ & $11$ \\  \hline

\multirow{10.5}{*}{$EIII$} & \multirow{10.5}{*}{$\begin{tikzpicture}[baseline=-1.8ex]\dynkin[ply=2,foldradius=2mm]{E}{oo***o}\end{tikzpicture}$} & \multirow{10.5}{*}{$41$} & $\begin{tikzpicture}[baseline=-0.7ex]\dynkin[ply=2,foldradius=2mm]{D}{*o*oo}\end{tikzpicture} \times \T$ & $28$ \\ 
&&& $\begin{tikzpicture}[baseline=-0.7ex]\dynkin[ply=2,foldradius=2mm]{D}{oo***}\end{tikzpicture} \times \T$ & $25$ \\
&&& $\begin{tikzpicture}[baseline=-1.8ex]\dynkin[ply=2,foldradius=2mm]{A}{oo*oo}\end{tikzpicture} \times \dynkin{A}{*}$ & $24$ \\
&&& $\dynkin{F}{II}$ & $24$ \\ 
&&& $\dynkin{C}{*o*o}$ & $22$ \\ 
&&& $\begin{tikzpicture}[baseline=-1.8ex]\dynkin[ply=2,foldradius=2mm]{A}{o***o}\end{tikzpicture} \times \dynkin{A}{o}$ & $21$ \\
&&& $(\dynkin[ply=2,foldradius = 2mm,rotate=90]{A}{oo})^2 \times \dynkin[ply=2,foldradius = 2mm,rotate=90]{A}{**}$ & $16$ \\
&&& $\begin{tikzpicture}[baseline=-0.7ex]\dynkin[ply=2,foldradius=2mm]{D}{*****}\end{tikzpicture} \times \T$ & $12$ \\
&&& $\dynkin{G}{2} \times \dynkin[ply=2,foldradius = 2mm,rotate=90]{A}{oo}$ & $11$ \\  \hline

\multirow{17.5}{*}{$EII$} & \multirow{17.5}{*}{$\begin{tikzpicture}[baseline=-1.8ex]\dynkin[ply=2,foldradius=2mm]{E}{oooooo}\end{tikzpicture}$} & \multirow{17.5}{*}{$46$} & $\dynkin{F}{I}$ & $32$ \\ 
&&& $\begin{tikzpicture}[baseline=-0.7ex]\dynkin[ply=2,foldradius=2mm]{D}{ooooo}\end{tikzpicture} \times \T$ & $30$ \\ 
&&& $\begin{tikzpicture}[baseline=-0.7ex]\dynkin[ply=2,foldradius=2mm]{D}{*o*oo}\end{tikzpicture} \times \T$ & $28$ \\
&&& $\begin{tikzpicture}[baseline=-1.8ex]\dynkin[ply=2,foldradius=2mm]{A}{ooooo}\end{tikzpicture} \times \dynkin{A}{o}$ & $26$ \\
&&& $\begin{tikzpicture}[baseline=-1.8ex]\dynkin[ply=2,foldradius=2mm]{A}{oo*oo}\end{tikzpicture} \times \dynkin{A}{*}$ & $24$ \\
&&& $\dynkin{C}{oooo}$ & $24$ \\ 
&&& $\dynkin{C}{*o**}$ & $20$ \\ 
&&& $\begin{tikzpicture}[baseline=0.35ex] 
\dynkin[name=upper]{A}{oo}
\node (current) at ($(upper root 1)+(-0,3.7mm)$) {};
\dynkin[at=(current),name=lower]{A}{oo}
\begin{scope}[on background layer]
\foreach \i in {1,2}%
{%
\draw[/Dynkin diagram/foldStyle]
($(upper root \i)$) -- ($(lower
root \i)$);%
}%
\end{scope}
\end{tikzpicture} \times \dynkin{A}{oo}$ & $19$ \\
&&& $(\dynkin[ply=2,foldradius = 2mm,rotate=90]{A}{oo})^3$ & $18$ \\
&&& $\dynkin{G}{I} \times \dynkin[ply=2,foldradius = 2mm,rotate=90]{A}{oo}$ & $16$ \\ 
&&& $\dynkin{G}{I} \times \dynkin[ply=2,foldradius = 2mm,rotate=90]{A}{**}$ & $14$ \\ 
&&& $(\dynkin[ply=2,foldradius = 2mm,rotate=90]{A}{**})^2 \times \dynkin[ply=2,foldradius = 2mm,rotate=90]{A}{oo}$ & $14$ \\
&&& $\begin{tikzpicture}[baseline=-1.8ex]\dynkin[ply=2,foldradius=2mm]{A}{*****}\end{tikzpicture} \times \dynkin{A}{*}$ & $12$ \\
&&& $\dynkin{G}{I}$ & $10$ \\ 
&&& $\dynkin{A}{oo}$ & $7$ \\ 
&&& $\dynkin[ply=2,foldradius = 2mm,rotate=90]{A}{oo}$ & $6$ \\ \hline

\multirow{14}{*}{$EVII$} & \multirow{14}{*}{$\dynkin{E}{o****oo}$} & \multirow{14}{*}{$66$} & $\begin{tikzpicture}[baseline=-1.8ex]\dynkin[ply=2,foldradius=2mm]{E}{oo***o}\end{tikzpicture} \times \T$ & $42$ \\*
&&& $\begin{tikzpicture}[scale=.175,baseline=0.3ex]
    \foreach \y in {0,...,3}
    \draw[thin,xshift=\y cm] (\y cm,0) ++(.3 cm, 0) -- +(14 mm,0);
    \draw[thin,fill=black] (0 cm,0) circle (3 mm);
    \draw[thin] (2 cm,0) circle (3 mm);
    \draw[thin,fill=black] (4 cm,0) circle (3 mm);
    \draw[thin] (6 cm,0) circle (3 mm);
    \draw[thin] (8 cm,0) circle (3 mm);
    \draw[thin,fill=black] (6 cm,2 cm) circle (3 mm);
    \draw[thin, fill=black] (6 cm, 3mm) -- +(0, 1.4 cm);
  \end{tikzpicture} ~\times \dynkin{A}{*}$ & $41$ \\
&&& $\begin{tikzpicture}[baseline=0.3ex]\dynkin{E}{o****o}\end{tikzpicture} \times \T$ & $38$ \\
&&& $\dynkin{A}{*o*o*o*}$ &$38$ \\
&&& $\begin{tikzpicture}[baseline=-1.8ex]\dynkin[ply=2,foldradius=2mm]{A}{oo***oo}\end{tikzpicture}$ &$34$ \\
&&& $\begin{tikzpicture}[scale=.175,baseline=0.3ex]
    \foreach \y in {0,...,3}
    \draw[thin,xshift=\y cm] (\y cm,0) ++(.3 cm, 0) -- +(14 mm,0);
    \draw[thin] (0 cm,0) circle (3 mm);
    \draw[thin] (2 cm,0) circle (3 mm);
    \draw[thin,fill=black] (4 cm,0) circle (3 mm);
    \draw[thin,fill=black] (6 cm,0) circle (3 mm);
    \draw[thin,fill=black] (8 cm,0) circle (3 mm);
    \draw[thin,fill=black] (6 cm,2 cm) circle (3 mm);
    \draw[thin] (6 cm, 3mm) -- +(0, 1.4 cm);
  \end{tikzpicture} ~\times \dynkin{A}{o}$ & $34$ \\
&&& $\begin{tikzpicture}[baseline=-1.8ex]\dynkin[ply=2,foldradius=2mm]{A}{ooooo}\end{tikzpicture}\times \dynkin[ply=2,foldradius = 2mm,rotate=90]{A}{**}$ & $27$ \\
&&& $\dynkin{F}{II} \times \dynkin{A}{o}$ & $27$ \\
&&& $\begin{tikzpicture}[baseline=-1.8ex]\dynkin[ply=2,foldradius=2mm]{A}{o***o}\end{tikzpicture}\times \dynkin[ply=2,foldradius = 2mm,rotate=90]{A}{oo}$ & $24$ \\
&&& $\dynkin{G}{2} \times \dynkin{C}{ooo}$ & $20$ \\ 
&&& $\begin{tikzpicture}[baseline=-1.8ex]\dynkin[ply=2,foldradius=2mm]{E}{******}\end{tikzpicture} \times \T$ & $14$ \\
&&& $\dynkin{F}{4} \times \dynkin{A}{o}$ & $14$ \\ \hline

\multirow{21}{*}{$EVI$} & \multirow{21}{*}{$\dynkin{E}{o*oo*o*}$} & \multirow{21}{*}{$74$} & $\begin{tikzpicture}[baseline=-1.8ex]\dynkin[ply=2,foldradius=2mm]{E}{oooooo}\end{tikzpicture} \times \T$ & $47$ \\ 
&&& $\begin{tikzpicture}[scale=.175,baseline=0.3ex]
    \foreach \y in {0,...,3}
    \draw[thin,xshift=\y cm] (\y cm,0) ++(.3 cm, 0) -- +(14 mm,0);
    \draw[thin,fill=black] (0 cm,0) circle (3 mm);
    \draw[thin] (2 cm,0) circle (3 mm);
    \draw[thin,fill=black] (4 cm,0) circle (3 mm);
    \draw[thin] (6 cm,0) circle (3 mm);
    \draw[thin] (8 cm,0) circle (3 mm);
    \draw[thin,fill=black] (6 cm,2 cm) circle (3 mm);
    \draw[thin, fill=black] (6 cm, 3mm) -- +(0, 1.4 cm);
  \end{tikzpicture} ~\times \dynkin{A}{o}$ & $42$ \\
&&& $\begin{tikzpicture}[scale=.175,baseline=0.3ex]
    \foreach \y in {0,...,3}
    \draw[thin,xshift=\y cm] (\y cm,0) ++(.3 cm, 0) -- +(14 mm,0);
    \draw[thin] (0 cm,0) circle (3 mm);
    \draw[thin] (2 cm,0) circle (3 mm);
    \draw[thin] (4 cm,0) circle (3 mm);
    \draw[thin] (6 cm,0) circle (3 mm);
    \draw[thin,fill=black] (8 cm,0) circle (3 mm);
    \draw[thin,fill=black] (6 cm,2 cm) circle (3 mm);
    \draw[thin] (6 cm, 3mm) -- +(0, 1.4 cm);
  \end{tikzpicture} ~\times \dynkin{A}{*}$ & $~42$ \topstrut \\
&&& $\begin{tikzpicture}[baseline=-1.8ex]\dynkin[ply=2,foldradius=2mm]{E}{oo***o}\end{tikzpicture} \times \T$ & $42$ \\
&&& $\begin{tikzpicture}[baseline=-1.8ex]\dynkin[ply=2,foldradius=2mm]{A}{ooooooo}\end{tikzpicture}$ & $39$ \\
&&& $\begin{tikzpicture}[baseline=-1.8ex]\dynkin[ply=2,foldradius=2mm]{A}{oo***oo}\end{tikzpicture}$ & $34$ \\
&&& $\dynkin{F}{I} \times \dynkin{A}{*}$ & $34$ \\
&&& $\dynkin{A}{*o*o*} \times \dynkin{A}{oo}$ & $29$ \\
&&& $\begin{tikzpicture}[baseline=-1.8ex]\dynkin[ply=2,foldradius=2mm]{A}{oo*oo}\end{tikzpicture}\times \dynkin[ply=2,foldradius = 2mm,rotate=90]{A}{oo}$ & $28$ \\
&&& $\begin{tikzpicture}[baseline=-1.8ex]\dynkin[ply=2,foldradius=2mm]{A}{oo*oo}\end{tikzpicture} \times \dynkin[ply=2,foldradius = 2mm,rotate=90]{A}{**}$ & $26$ \\
&&& $\dynkin{F}{II} \times \dynkin{A}{*}$ & $26$ \\
&&& $\begin{tikzpicture}[scale=.175,baseline=0.3ex]
    \foreach \y in {0,...,3}
    \draw[thin,xshift=\y cm] (\y cm,0) ++(.3 cm, 0) -- +(14 mm,0);
    \draw[thin,fill=black] (0 cm,0) circle (3 mm);
    \draw[thin,fill=black] (2 cm,0) circle (3 mm);
    \draw[thin,fill=black] (4 cm,0) circle (3 mm);
    \draw[thin,fill=black] (6 cm,0) circle (3 mm);
    \draw[thin,fill=black] (8 cm,0) circle (3 mm);
    \draw[thin,fill=black] (6 cm,2 cm) circle (3 mm);
    \draw[thin,fill=black] (6 cm, 3mm) -- +(0, 1.4 cm);
  \end{tikzpicture} ~\times \dynkin{A}{*}$ & $24$ \\
&&& $\dynkin{G}{oo} \times \dynkin{C}{*o*}$ & $23$ \\ 
&&& $\dynkin{G}{2} \times \dynkin{C}{*o*}$ & $18$ \\ 
&&& $\dynkin{G}{oo} \times \dynkin{C}{***}$ & $18$ \\ 
&&& $\begin{tikzpicture}[baseline=-1.8ex]\dynkin[ply=2,foldradius=2mm]{A}{*****}\end{tikzpicture}\times \dynkin[ply=2,foldradius = 2mm,rotate=90]{A}{oo}$ & $16$ \\
&&& $\dynkin{G}{oo} \times \dynkin{A}{*}$ & $12$ \\
&&& $\dynkin[ply=2,foldradius = 2mm,rotate=90]{A}{oo}$ & $6$ \\ 
&&& $\dynkin{A}{o} \times \dynkin{A}{*}$ & $5$ \\   \hline

\multirow{21.5}{*}{$EIX$} & \multirow{21.5}{*}{$\dynkin{E}{o****ooo}$} & \multirow{21.5}{*}{$125$} & $\begin{tikzpicture}[baseline=0.3ex]\dynkin{E}{o*oo*o*}\end{tikzpicture} \times \dynkin{A}{1}$ & $76$ \\
&&& $\begin{tikzpicture}[baseline=0.3ex]\dynkin{E}{o****oo}\end{tikzpicture} \times \dynkin{A}{o}$ & $69$ \\
&&& $\begin{tikzpicture}[scale=.175,baseline=0.3ex]
    \foreach \y in {0,...,5}
    \draw[thin,xshift=\y cm] (\y cm,0) ++(.3 cm, 0) -- +(14 mm,0);
    \draw[thin,fill=black] (0 cm,0) circle (3 mm);
    \draw[thin] (2 cm,0) circle (3 mm);
    \draw[thin,fill=black] (4 cm,0) circle (3 mm);
    \draw[thin] (6 cm,0) circle (3 mm);
    \draw[thin,fill=black] (8 cm,0) circle (3 mm);
    \draw[thin] (10 cm,0) circle (3 mm);
    \draw[thin] (12 cm,0) circle (3 mm);
    \draw[thin,fill=black] (10 cm,2 cm) circle (3 mm);
    \draw[thin,fill=black] (10 cm, 3mm) -- +(0, 1.4 cm);
  \end{tikzpicture}$ & $68$ \\
&&& $\begin{tikzpicture}[scale=.175,baseline=0.3ex]
    \foreach \y in {0,...,5}
    \draw[thin,xshift=\y cm] (\y cm,0) ++(.3 cm, 0) -- +(14 mm,0);
    \draw[thin] (0 cm,0) circle (3 mm);
    \draw[thin] (2 cm,0) circle (3 mm);
    \draw[thin] (4 cm,0) circle (3 mm);
    \draw[thin] (6 cm,0) circle (3 mm);
    \draw[thin,fill=black] (8 cm,0) circle (3 mm);
    \draw[thin,fill=black] (10 cm,0) circle (3 mm);
    \draw[thin,fill=black] (12 cm,0) circle (3 mm);
    \draw[thin,fill=black] (10 cm,2 cm) circle (3 mm);
    \draw[thin] (10 cm, 3mm) -- +(0, 1.4 cm);
  \end{tikzpicture}$ & $~61$ \topstrut \\
&&& $\begin{tikzpicture}[baseline=-1.8ex]\dynkin[ply=2,foldradius=2mm]{E}{oooooo}\end{tikzpicture} \times \dynkin[ply=2,foldradius = 2mm,rotate=90]{A}{**}$ & $50$ \\
&&& $\begin{tikzpicture}[baseline=-1.8ex]\dynkin[ply=2,foldradius=2mm]{E}{oo***o}\end{tikzpicture} \times \dynkin[ply=2,foldradius = 2mm,rotate=90]{A}{oo}$ & $47$ \\
&&& $\begin{tikzpicture}[baseline=-1.8ex]\dynkin[ply=2,foldradius=2mm]{A}{ooo**ooo}\end{tikzpicture}$ & $46$ \\
&&& $\begin{tikzpicture}[baseline=0.3ex]\dynkin{E}{o****o}\end{tikzpicture}\times \dynkin{A}{oo}$ & $44$ \\
&&& $\begin{tikzpicture}[baseline=-1.8ex]\dynkin[ply=2,foldradius=2mm]{A}{oo****oo}\end{tikzpicture}$ & $40$ \\
&&& $\dynkin{G}{2} \times \dynkin{F}{I}$ & $37$ \\ 
&&& $\dynkin{G}{I} \times \dynkin{F}{II}$ & $34$ \\ 
&&& $\begin{tikzpicture}[baseline=-1.8ex]\dynkin[ply=2,foldradius=2mm]{A}{oooo}\end{tikzpicture} \times  \begin{tikzpicture}[baseline=-1.8ex]\dynkin[ply=2,foldradius=2mm]{A}{o**o}\end{tikzpicture}$ & $30$ \\
&&& $\begin{tikzpicture}[baseline=-1.8ex]\dynkin[ply=2,foldradius=2mm]{A}{oooo}\end{tikzpicture} \times  \begin{tikzpicture}[baseline=-1.8ex]\dynkin[ply=2,foldradius=2mm]{A}{****}\end{tikzpicture}$ & $24$ \\
&&& $\dynkin{G}{I} \times \dynkin{F}{4}$ & $21$ \\ 
&&& $\begin{tikzpicture}[baseline=0.4ex] 
\dynkin[name=upper]{G}{oo}
\node (current) at ($(upper root 1)+(-0,3.7mm)$) {};
\dynkin[at=(current),name=lower]{G}{oo}
\begin{scope}[on background layer]
\foreach \i in {1,2}%
{%
\draw[/Dynkin diagram/foldStyle]
($(upper root \i)$) -- ($(lower
root \i)$);%
}%
\end{scope}
\end{tikzpicture}\times \dynkin{A}{1}$ & $20$ \\
&&& $\begin{tikzpicture}[baseline=0.3ex]\dynkin{E}{*******}\end{tikzpicture} \times \dynkin{A}{1}$ & $19$ \\
&&& $\begin{tikzpicture}[baseline=-1.8ex]\dynkin[ply=2,foldradius=2mm]{E}{******}\end{tikzpicture} \times \dynkin[ply=2,foldradius = 2mm,rotate=90]{A}{oo}$ & $19$ \\
&&& $\dynkin{A}{oo} \times \dynkin{A}{1}$ & $9$ \\ \hline

\end{longtable}
 \restoregeometry

\section{Proof of Theorem \ref{depth}}

\begin{proof} Let $X$ be a simple complex algebraic group and let $\mathcal{X}:= \big\{ X=X_0 > X_1 >  ... > X_k=H \big\}$ be an unrefinable chain of connected reductive complex subgroups of $X$ (where $H$ is not necessarily the trivial group). For every $i$ let $(X_i)_s$ be a split form of $X_i$. Recall from Proposition \ref{splitform} that such an $(X_i)_s$ always exists and is unique up to conjugacy in $X_i$. If $(X_i)_s$ is contained in some conjugate of $(X_{i+1})_s$ for every $i$ then we say that the chain $\mathcal{X}$ \textit{splits}. That is, after adjusting by an appropriate set of conjugates, there exists a chain $\mathcal{X}_s := \big\{(X_0)_s > (X_1)_s > (X_2)_s > ... > (X_k)_s \big\}$ of split real subgroups. Observe that $\mathcal{X}_s$ is unrefinable by Lemma \ref{maxconn}. 

\vspace{2mm}\noindent Now let $\mathcal{X}(X) = X > ... > H$ be the (unique) subchain of one of the following complex chains that satisfies $H=A_1$: $$\begin{cases} A_{2k-1} > C_k > A_1 & (k \geq 2) \\ A_{2k} > B_k > A_1 & (k \geq 4 \textnormal{ or } k=2) \\  A_6 > B_3 > G_2 > A_1 \\ D_{4} > A_2 > A_1 \\ D_{k} > B_{k-1} > A_1 & (k \geq 5) \\ E_6 > F_4> A_1 \\ E_7 > A_1  \\ E_8 > A_1  \end{cases}$$ Observe that $\mathcal{X}(X)$ has length $\lambda_{\C} (X)-3$ by Theorem \ref{BLSTHM3}. 

\begin{Lemma}\label{descends} Let $X$ be a simple complex group other than $D_4$. Then the chain $\mathcal{X}(X)$ splits.

\begin{proof} If $X=G_2$, $E_7$ or $E_8$ then $\mathcal{X}(X)$ splits by Tables $6$, $38-39$ and $55-57$ of \cite{K} respectively. The chain $\mathcal{X}(E_6)=E_6 > F_4 > A_1$ splits by Tables $12$ and $29$ of \cite{K}. The embedding $\mathrm{SO}_{2k}(\C) > \mathrm{SO}_{2k-1}(\C)$ splits into $\mathrm{SO}(k,k) > \mathrm{SO}(k,k-1)$ by definition of the indefinite orthogonal group. That is, the chain $D_k > B_{k-1}$ splits.

\vspace{2mm}\noindent Now let $R=(r_{ij})$ be the $n \times n$ matrix with entries given by $r_{ij}= i$ if $j=i+1$ and $r_{ij}= 0$ otherwise. Similarly, let $S=(s_{ij})$ be the $n \times n$ matrix with entries given by $s_{ij}=n-j$ if $j=i-1$ and $s_{ij}=0$ otherwise. In addition, let $P=(p_{ij})$ be a $n \times n$ antidiagonal matrix with non-zero entries that satisfy $ip_{i(n-i+1)}+(n-i)p_{(i+1)(n-i)}=0$ for all $1 \leq i \leq n-1$. For example, if $n=3$ then $P$ is a scalar multiple of $\begin{psmallmatrix} 0&0&2 \\0& \hspace{-0.7mm}-\hspace{-0.2mm}1 & 0 \\ 2&0&0 \end{psmallmatrix}$. Note that $P$ is symmetric if $n$ is odd and skew-symmetric if $n$ is even. 

\vspace{2mm}\noindent Let $V=\C^n$ be equipped with a non-degenerate bilinear form $\mathcal{P}$ that corresponds to the matrix $P$. Henceforth let $X=\SL_n(\C)$ and let $Y=\big\{A \in X \hspace{0.4mm} \big| \hspace{0.4mm} A^{\top}PA=P \big\}$ be the subgroup of $X$ that preserves $\mathcal{P}$. Let $H = A_1$ be the irreducibly embedded complex subgroup of $X$ that is generated by the matrices $\Exp(tR)$ and $\Exp(tS)$ for all $t \in \C$. We check that $R^{\top}P + PR=0=S^{\top}P + PS$ and hence $H<Y$ by Lemma $11.2.2$ of \cite{C1}. Let $\sigma$ be the antiholomorphic involution of $X$ that sends $A \mapsto \overline{A}$. Then $\sigma$ stabilises the chain $X > Y > H$ with fixed points $$\begin{cases} \SL_{2k+1}(\R) > \mathrm{SO}(k+1,k) > \PSL_2(\R) & \textnormal { if } n=2k+1\geq 3 \\ \SL_{2k}(\R) > \Sp_{2k}(\R) > \SL_2(\R) & \textnormal { if } n=2k \geq 4  \end{cases}$$ That is, the chains $A_{2k} > B_k > A_1$ and $A_{2k-1} > C_k > A_1 $ both split. If $n=7$ then there exists a $\sigma$-stable copy of $G_2$ in $Y$ that contains $H$. Hence the chain $A_{6} > B_3 > G_2 > A_1$ splits since the only non-compact real form of $G_2$ is the split form.
\end{proof}
\end{Lemma}

\begin{Lemma}\label{quasisplit} Let $G$ be a quasisplit simple real algebraic group. Then $\lambda(G) \leq \lambda_{\C} \big(G(\C)\big)-1$. 
\begin{proof} We first consider the case where $G=X_{\R}$ for some simple complex group $X$. Recall that $G(\C) \cong X^2$ and so $\lambda_{\C} \big(G(\C)\big) = \lambda_{\C}(X)+1$ by Corollary \ref{reductiondepthcor}. The maximal compact subgroup of $G$ is isomorphic to the compact form $X_c$ of $X$. So $\lambda(G) \leq \lambda(X_c) +1 =\lambda_{\C} (X) = \lambda_{\C} \big(G(\C)\big)-1$ by Proposition \ref{maxcompact} and Theorem \ref{BLSTHM4}.

\vspace{2mm}\noindent Henceforth we can assume that $G(\C)$ is a simple complex group. The maximal compact subgroup of $(A_1)_s$ is a compact torus $\T$. So the real chain $(A_1)_s > \T > 1$ is unrefinable by Proposition \ref{maxcompact}. If $G$ is split and not of type $D_4$ then there exists an unrefinable chain $$G > ... > (A_1)_s > \T > 1$$ of length $\lambda_{\C} \big(G(\C)\big)-1$ by Lemma \ref{descends}. If $G$ is split and of type $D_4$ then the chain $$ \mathrm{SO}(4,4) > \SU(2,1) > \mathrm{SO}(2,1) > \T >1 $$ is unrefinable by Table $19$ of \cite{K} and Lemma \ref{maxconn}.

\vspace{2mm}\noindent It remains to consider the cases where $G$ is quasisplit but not split. If $G$ is of type $A_{2k}$ for $k \geq 1$ (resp. $A_{2k-1}$ for $k \geq 2$, $D_k$ for $k \geq 5$, $D_4$ or $E_6$) then by Lemmas \ref{maxconn} and \ref{descends} it suffices to check that $G$ contains a split form of $B_k$ (resp. $C_{k}$, $B_{k-1}$, $A_2$ or $F_4$).

\vspace{2mm}\noindent For any $k$ the embeddings $\SU(k+1,k) > \mathrm{SO}(k+1,k)$ and $\mathrm{SO}(k+1,k-1) > \mathrm{SO}(k,k-1)$ follow immediately from their definitions. The embeddings $\mathrm{SO}(5,3) > \SL_3(\R)$ and $EII > FI$ are given in Tables $12$ and $19$ of \cite{K} respectively.

\vspace{2mm}\noindent Finally, let $Y=\big\{A \in \SL_{2k}(\C) \hspace{0.4mm} \big| \hspace{0.4mm} A^{\top}QA = Q \big\} \cong \Sp_{2k}(\C)$ for $k \geq 2$ where $Q :=\begin{psmallmatrix} 0&I_k \\ \hspace{-0.7mm}-\hspace{-0.2mm}I_k&0\end{psmallmatrix}$ and $I_k$ denotes the $k$ by $k$ identity matrix. Let $s := \diag\!\big(i^{(k)},(-i)^{(k)}\big) \in Y$. The antiholomorphic involution $A \mapsto s\big({\overline{A}}^{\top}\big)^{-1}s^{-1}$ of $\SL_{2k}(\C)$ stabilises $Y$ with fixed points $\SU(k,k) > \Sp_{2k}(\R)$.
\end{proof}
\end{Lemma}

\noindent The following Lemma is similar to Lemma $2.3$ of \cite{BLS}.

 \begin{Lemma}\label{lowdim} Let $G$ be a real algebraic group. Then $\lambda(G)=1$ if and only if $\dim(G)=1$ and $\lambda(G)=2$ if and only if $\dim(G)=2$ or $G(\C)=A_1$. 
\begin{proof} The case where $\lambda(G)=1$ is obvious. So assume that $\lambda(G)=2$. If $G$ is soluble then $\dim(G)=2$ by Lemma $2.2$ of \cite{BLS}. If $G$ is insoluble then $G' \neq 1$ and $\lambda(G') \geq 2$. Applying Lemma \ref{lemma2.5} to the decomposition $G=R(G) \cdot G'$ implies that $\lambda(G') = 2$ and $R(G) =1$. 

\vspace{2mm}\noindent If $G$ is compact then $r(G)=1$ by Theorem \ref{BLSTHM4} and Corollary \ref{reductiondepth}. So we can assume that $G$ is not compact. Let $M$ be a maximal connected subgroup of $G$ of dimension $1$. If $M$ is unipotent then $M$ is strictly contained in a parabolic subgroup of $G$, a contradiction. If $M$ is a torus then $M$ is contained in a maximal torus $T$ of $G$ with $\dim(T)=r(G)$ and hence $r(G)=1$.

\vspace{2mm}\noindent There are two real forms of $A_1$ up to isogeny, a split form $(A_1)_s$ and a compact form $(A_1)_c$. The chain $(A_1)_s > \T >1$ is unrefinable by Proposition \ref{maxcompact} and $\lambda \big((A_1)_c\big)=2$ by Theorem \ref{BLSTHM4}.
\end{proof}
\end{Lemma}

\noindent Now let $U_k$ denote a unipotent group of dimension $k$ and similarly let $T^k$ denote a $k$-dimensional torus.

\begin{Lemma}\label{depthis3} Let $G$ be a non-compact semisimple real algebraic group. Then $\lambda(G)=3$ if and only if $G$ is quasisplit and is of type $(A_1)^2$, $A_2$, $C_r$ ($r \geq 2$), $B_r$ ($r \geq 4$), $G_2$, $F_4$, $E_7$ or $E_8$.
\begin{proof} By Lemma \ref{lowdim} we can assume that $r(G)>1$. Let $\lambda(G)=3$ and let $M$ be a maximal connected subgroup of $G$ such that $\lambda(M)=2$. Assume (for a contradiction) that $\dim(M)=2$. If $M=\C^{\times}$, $(\R^\times)^2$ or if $M$ is unipotent then $M$ is strictly contained in a minimal parabolic subgroup of $G$. If $M=\T^2$ then $M$ is strictly contained in a maximal compact subgroup of $G$. If $M=U_1T^1$ then $M<N_G(U_1)\leq P$ for some parabolic subgroup $P$ of $G$ (Borel-Tits, Theorem $2.5$ of \cite{BT}). All of these cases contradict the maximality of $M$. Hence $M$ is of type $A_1$ by Lemma \ref{lowdim}.

\vspace{2mm}\noindent If $G$ is not simple then $G$ must be of type $(A_1)^2$ where $M$ is diagonally embedded in $G$ by Corollary \ref{reductiondepth} and Lemma \ref{lowdim}. Otherwise, $G$ is simple and so either $G=X_{\R}$ for some simple complex group $X$ or $G(\C)$ is simple. If $G=X_{\R}$ then $X=A_1$ by Proposition \ref{simplenotsimple} and so again $G$ is of type $(A_1)^2$. If $G(\C)$ is simple then $M(\C)$ is maximal connected in $G(\C)$ by Lemmas $2.4$ and $2.7$ of \cite{K} (since $M$ is of type $A_1$). Hence either $G(\C)=(A_1)^2$ or $G(\C)$ is simple and contains a maximal connected copy of $A_1$. That is, $G(\C)$ is one of $(A_1)^2$, $A_2$, $C_r$ ($r \geq 2$), $B_r$ ($r \geq 4$), $G_2$, $F_4$, $E_7$ or $E_8$.

\vspace{2mm}\noindent Without loss of generality we may assume that $G(\C)$ is adjoint. Observe that $C_{G(\C)}\big(M(\C)\big)$ is finite since $M(\C)$ is maximal connected and $Z\big(M(\C)\big)$ is finite. Hence all elements of $C_{G(\C)}\big(M(\C)\big)$ are semisimple. If $z \in C_{G(\C)}\big(M(\C)\big)$ is non-trivial then $M(\C)= {C_{G(\C)}(z)}^{\circ}<G(\C)$ by maximality. But ${C_{G(\C)}(z)}^{\circ}$ contains some maximal torus of $G(\C)$ by Theorem $14.2$ of \cite{MT}, which is a contradiction. Hence $ C_{G(\C)}\big(M(\C)\big)=1$ and in particular $M(\C) \cong \PSL_2(\C)$.

\vspace{2mm}\noindent Now let $G$ be an inner real form of $(A_1)^2$, $A_2$, $C_r$ ($r \geq 2$), $B_r$ ($r \geq 4$), $G_2$, $F_4$, $E_7$ or $E_8$. For a given $M(\C)$ and $G(\C)$, it follows from Theorem \ref{Cartan} that conjugacy classes of embeddings of real forms $M <G$ are in bijection with conjugacy classes of holomorphic inner involutions $\theta$ of $G(\C)$ that stabilise $M(\C)$. We can ignore the case where $\theta$ is trivial as this corresponds to the case where $M$ and $G$ are both compact. Since $ C_{G(\C)}\big(M(\C)\big)=1$, the only remaining possibility for $\theta$ is conjugation by an element $s \in M(\C)$ which has order $2$. But all elements of order $2$ in $M(\C) \cong \PSL_2(\C)$ are conjugate. So $\theta$ corresponds to the case where $M$ and $G$ are both split by Lemma \ref{descends}.

\vspace{2mm}\noindent It remains to consider the cases where $G$ is an outer real form. Then $\Out\!\big(G(\C)\big)$ is non-trivial and so $G$ can only be of type $A_2$ or $(A_1)^2$. The split form $(A_2)_s$ is the unique outer real form of $A_2$ and the unique outer real form $(A_1)_{\R}$ of $(A_1)^2$ is quasisplit. Both $(A_2)_s$ and $(A_1)_{\R}$ contain a maximal connected subgroup of type $A_1$ by Lemmas \ref{descends} and \ref{quasisplit} respectively.

\vspace{2mm}\noindent The converse follows from Theorem \ref{BLSTHM3} and Lemmas \ref{quasisplit} and \ref{lowdim}.
\end{proof} 
\end{Lemma}

\begin{Lemma}\label{depthis4} Let $G$ be quasisplit and of type $A_r ~(r \geq 3, r \neq 6)$, $D_r ~(r \geq 4)$, $B_3$ or $E_6$, or let $G=\SL_k(\mathbb{H}) ~(k>1)$, $\mathrm{SO}(2k\hspace{-0.5mm}+\hspace{-0.5mm}1,1) ~(k >3)$, $FII$, $EIV$ or $EVI$. Then $\lambda(G)=4$.
\begin{proof} It suffices to show that $\lambda(G) \leq 4$ since $G$ is a non-compact simple real group that is not of type $A_1$ and is not one of the groups listed in Lemma \ref{depthis3}. If $G$ is quasisplit then $\lambda(G) \leq 4$ by Lemma \ref{quasisplit} and Theorem \ref{BLSTHM3}.

\vspace{2mm}\noindent If $G=EVI$ then by Table \ref{exceptcheck} there exists a maximal connected $M \cong \PSU(2,1)$ in $G$. Then $\lambda(M)=3$ by Lemma \ref{depthis3} since $M$ is quasisplit and of type $A_2$.

\vspace{2mm}\noindent For all remaining cases, let $K$ be a maximal compact subgroup of $G$. If $G=\SL_k(\mathbb{H}) $ for $k>1$ (resp. $\mathrm{SO}(2k+1,1)$ for $k > 3$, $FII$, $EIV$) then ${K}^{\circ}=(C_k)_c$ (resp. $(B_{k})_c$, $(B_4)_c$, $(F_4)_c$). For each case $\lambda({K}^{\circ})=3$ by Theorems \ref{BLSTHM3} and \ref{BLSTHM4} and hence $\lambda(G) \leq 4$ by Proposition \ref{maxcompact}.
\end{proof}
\end{Lemma}

\begin{Lemma}\label{depthis5} If $G=EIII$, $EVII$ or $EIX$ then $\lambda(G)=5$. If $G$ is non-compact of type $A_6$ then $\lambda(G) \geq 5$ with equality if $G$ is quasisplit.
\begin{proof} We have $\lambda(G) \geq 4$ by Lemmas \ref{lowdim} and \ref{depthis3}. Let $M$ be a maximal connected subgroup of $G$.

\vspace{2mm}\noindent If $G$ is of type $A_6$ then by Table $3$ of \cite{K} either $M$ is parabolic, $M$ is of type $A_kA_{6-k-1}\T$ or $M$ is simple and irreducibly embedded in $G$ (in particular, $M$ is of type $A_1$, $G_2$ or $B_3$). If $M$ is simple and irreducibly embedded in $G$ then $M(\C)$ is maximal connected in $G(\C)$ by Lemma $2.4$ of \cite{K}, and so $M$ must be of type $B_3$.

\vspace{2mm}\noindent We first consider the case where $M$ is a maximal parabolic subgroup of $G$. Let $M=R_u(M) \rtimes L$ where $L$ is a Levi subgroup of $M$. Then $\lambda(M) \geq \lambda(L) +1 \geq 4$ by Corollary \ref{reductiondepth} and Lemma \ref{lowdim} since $\dim(L) >2$ and $L$ is not of type $A_1$.

\vspace{2mm}\noindent Now assume that $M$ is reductive and $\lambda(M)=3$. Then one of the following possibilities occurs. If $M$ is not semisimple then $\lambda(M) \geq \lambda(M')+1$ by Lemma \ref{lemma2.5} (since $R(M) \lhd M$) and so either $M'$ is trivial or $M'$ is of type $A_1$ by Lemma \ref{lowdim}. If $M$ is semisimple but not simple then $M=M_1M_2$ where $M_1$ and $M_2$ are both simple groups of type $A_1$ (by Corollary \ref{reductiondepth} and Lemma \ref{lowdim}). If $M$ is simple then $M$ is compact or quasisplit by Lemma \ref{depthis3}. If $G=EIII$, $EVII$ or $EIX$ then it is easy to see from Table \ref{exceptcheck} that none of these possibilities can occur. If $G$ is of type $A_6$ then $M$ must be of type $B_3$ but then $\lambda(M) \neq 3$ by Lemma \ref{depthis3}. We have a contradiction.

\vspace{2mm}\noindent It remains to show that $\lambda(G) \leq 5$ if $G=EIII$, $EVII$, $EIX$, $\SL_7(\R)$ or $\SU(4,3)$. If $G=\SL_7(\R)$ or $\SU(4,3)$ this follows from Theorem \ref{BLSTHM3} and Lemma \ref{quasisplit}. For each remaining case we use Table \ref{exceptcheck} and Lemma \ref{depthis4} to find a maximal connected subgroup $M$ of $G$ with $\lambda(M)=4$. If $G=EIII$ then $M=FII$ is such a subgroup and if $G=EVII$ then we take $M=\SL_4(\mathbb{H})$. If $G=EIX$ then let $M=\PSL_3(\R) \times \PSU(2)$ and observe that the chain $$G >M=\PSL_3(\R) \times \PSU(2) > \big(\!\PSU(2)\big)^2  > \PSU(2) > \T >1$$ is unrefinable by Proposition \ref{maxcompact} and Corollary \ref{reductiondepthcor}. 
\end{proof}
\end{Lemma}

\noindent By combining Theorem \ref{BLSTHM3} with Lemmas \ref{lowdim}, \ref{depthis3}, \ref{depthis4} and \ref{depthis5}, for $G(\C)$ a simple complex group, we have shown that $\lambda(G) \geq \lambda_{\C} \big(G(\C)\big)-1$ with equality if $G$ is quasisplit and proved parts (i) and (ii) of Theorem \ref{depth}. The following lemma takes care of the case where $G(\C)$ is not simple.

\begin{Lemma}\label{weirdones} Let $G=X_{\R}$ for $X$ a simple complex group. Then $\lambda(G)=\lambda_{\C} \big(G(\C)\big)-1=\lambda_{\C}(X)$.
\begin{proof} The case where $X=A_1$ has been done in Lemma \ref{depthis3}. Recall that $G(\C) \cong X^2$ and so $\lambda_{\C} \big(G(\C)\big) = \lambda_{\C}(X)+1$ by Corollary \ref{reductiondepthcor}. The upper bound $\lambda(G)  \leq \lambda_{\C} (X)$ was shown in Lemma \ref{quasisplit} since $G$ is quasisplit. It remains to show that $\lambda(G) \geq \lambda_{\C} (X)$ for $X$ a simple complex group with $r(X)>1$. Let $M$ be a maximal connected subgroup of $G$. Then one of the following three possibilities occurs by Proposition \ref{simplenotsimple}. Either $M$ is parabolic, $M$ is a real form of $X$ or $M=H_{\R}$ for some reductive maximal connected complex subgroup $H$ of $X$.

\vspace{2mm}\noindent Let $M$ be a maximal parabolic subgroup of $G$ and recall that $r(X)>1$. Then $M/R_u(M) =: L \cong L' \times \C^{\times}$ where $\C^{\times}$ is a real torus of dimension $2$ and $r(L') > 1$. Hence $\lambda(M) \geq \lambda(L)+1 \geq \lambda(L')+3 \geq 6$ by Lemma \ref{lemma2.5}, Corollary \ref{reductiondepth} and Lemma \ref{lowdim} since $L'$ is not of type $A_1$. So $\lambda(M) \geq 6 \geq \lambda_{\C}(X)$ by Theorem \ref{BLSTHM3}.

\vspace{2mm}\noindent If $M$ is a real form of $X$ then $\lambda(M) \geq \lambda_{\C}(X)-1$ by combining Lemmas \ref{lowdim}, \ref{depthis3} and \ref{depthis5}.

\vspace{2mm}\noindent Finally, let $M=H_{\R}$ for $H<X$ as above. Observe that $H'$ is non-trivial since $r(X)>1$. Hence $\lambda(M) \geq 3$ with equality if and only if $H=A_1$ by Lemmas \ref{lemma2.5}, \ref{lowdim} and \ref{depthis3}. So it suffices to consider the cases where $\lambda_{\C}(X) \geq 5$. If $X=A_r ~(r \geq 3, r \neq 6)$, $D_r ~(r \geq 4)$, $B_3$ or $E_6$ then $\lambda_{\C}(X)=5$ by Theorem \ref{BLSTHM3} and $\lambda(M) \geq 4$ since $X$ does not contain a maximal connected copy of $A_1$. By Theorem \ref{BLSTHM3} the only remaining case is $X=A_6$, which satisfies $\lambda_{\C}(X)=6$. Then $H = B_3$ by $\S 18$ of \cite{MT}, and $\lambda\big((B_3)_{\R}\big)\geq 5$ by the preceding arguments.
\end{proof}
\end{Lemma}

\vspace{1mm}\noindent It remains to show the upper bounds for classical $G$ given in part (iii) of Theorem \ref{depth}. We use Theorem \ref{Komclass1} and Proposition \ref{maxcompact} to construct unrefinable chains for classical $G$.

\vspace{2mm}\noindent Firstly, let $G=\mathrm{SO}^*(2k)$ for $k \geq 4$. By Theorem \ref{Komclass1} there exists a maximal connected $M \cong \mathrm{SO}_k(\C)_{\R}$ in $G$. The maximal compact subgroup of $M$ is isomorphic to $\mathrm{SO}(k)$ and so $\lambda(G) \leq \lambda(M)+1\leq  \lambda\big(\mathrm{SO}(k)\big) + 2 = 6-\zeta_k$ by Theorem \ref{BLSTHM4}.

\vspace{2mm}\noindent Next let $G=\Sp(p,q)$ for $p \geq q >0$. If $q>1$ then the chain $$G>\Sp(p)\times\Sp(q) > \Sp(p) \times \Sp(1) > \big(\!\Sp(1)\big)^2 > \Sp(1) > \T >1$$ is unrefinable and so $\lambda(G)\leq 6$. Now if $p>q=1$ then $$G>\Sp(p)\times\Sp(1) > \big(\!\Sp(1)\big)^2 > \Sp(1) > \T >1$$ is unrefinable and so $\lambda(G)\leq 5$. If $p=q>1$ then $$G>\big(\!\Sp(p)\big)^2 > \Sp(p) > \SU(2) > \T >1$$ is unrefinable and so again $\lambda(G)\leq 5$. Finally, if $p=q=1$ then $$G>\big(\!\Sp(1)\big)^2 > \Sp(1) >\T >1$$ is unrefinable and hence $\lambda(G) = 4$ by Lemmas \ref{lowdim} and \ref{depthis3}.

\vspace{2mm}\noindent Now let $G=\SU(p,q)$ for $p \geq q >0$. Then $G$ contains a copy of $\mathrm{SO}(p,q)$ which is maximal connected by Lemma \ref{maxconn}. Hence $\lambda(G) \leq \lambda\big(\mathrm{SO}(p,q) \big)+1$. 

\vspace{2mm}\noindent It remains to consider the most complicated case, let $G=\mathrm{SO}(p,q)$ for $p \geq q >0$. For any choice of integers satisfying $p_1+p_2=p$ and $q_1+q_2=q$, recall from Table \ref{classicalcheck} that (the connected component of) $\mathrm{SO}(p_1,q_1) \times \mathrm{SO}(p_2,q_2)$ is a maximal connected subgroup of $G$.

\vspace{2mm}\noindent If $p-q \leq 2$ then $G$ is quasisplit and so $\lambda(G) \leq 4$ by Lemma \ref{quasisplit}. If $p-q =3$ then $M=\mathrm{SO}(p-1,q)$ is a quasisplit maximal connected subgroup of $G$ and so $\lambda(G)\leq 5$. Similarly, if $p-q =4$ then $M=\mathrm{SO}(p-1,q)$ is a maximal connected subgroup of $G$ with $\lambda(M) \leq 5$ and so $\lambda(G)\leq 6$. Note that if $p=12$ and $q=7$ then $\lambda(G) \leq \lambda\big(\mathrm{SO}(11,7)\big) +1\leq 7$.

\vspace{2mm}\noindent If $q=0$ then $G$ is compact and so $\lambda(G) \leq 4$ by Theorem \ref{BLSTHM4}. If $q =1$ then $M=\mathrm{SO}(p)$ is a compact maximal connected subgroup of $G$ and so $\lambda(G)\leq 5$. Similarly, if $q =2$ then $M=\mathrm{SO}(p,1)$ is a maximal connected subgroup of $G$ with $\lambda(G) \leq 5$ and so $\lambda(G)\leq 6$. 

\vspace{2mm}\noindent So henceforth we can assume that $p-q >4 $ and $q>2$. In particular, $p \neq 7$.

\vspace{2mm}\noindent If $p$ is odd and $q =7$ then $$G> \mathrm{SO}(p) \times \mathrm{SO}(7) > (A_1)_c \times \mathrm{SO}(7) > (A_1)_c \times (G_2)_c > (A_1)_c^2 > (A_1)_c > \T >1$$ is unrefinable and so $\lambda(G) \leq 7$. If $p$ and $q$ are odd and $q \neq 7$ then the chain $$G > \mathrm{SO}(p) \times \mathrm{SO}(q) > (A_1)_c \times \mathrm{SO}(q) > (A_1)_c^2 > (A_1)_c > \T >1$$ is unrefinable and so $\lambda(G) \leq 6$. So if $p$ is even and $q \neq 7$ is odd then $\lambda(G) \leq \lambda\big(\mathrm{SO}(p-1,q)\big) +1\leq 7$. Similarly, if $p$ is odd and $q \neq 8$ is even then $\lambda(G) \leq \lambda\big(\mathrm{SO}(p,q-1)\big) +1\leq 7$. If $p$ is odd and $q =8$ then $$G> \mathrm{SO}(p) \times \mathrm{SO}(8) > (A_1)_c \times \mathrm{SO}(8) > (A_1)_c \times (A_2)_c > (A_1)_c^2 > (A_1)_c > \T >1$$ is unrefinable and so again $\lambda(G) \leq 7$.

\vspace{2mm}\noindent If $p$ and $q$ are even and $p-q \neq 8$ then the chain $$G > \mathrm{SO}(p-q - 1) \times \mathrm{SO}(q + 1,q) > (A_1)_c \times \mathrm{SO}(q + 1,q) > (A_1)_c \times (A_1)_s > (A_1)_c \times \T > \T^2 > \T >1$$ is unrefinable and so $\lambda(G) \leq 7$. Similarly, if $p$ and $q$ are even and $p-q = 8$ but $q \neq 4$ then $$G> \mathrm{SO}(p-q + 1) \times \mathrm{SO}(q - 1,q) > (A_1)_c \times \mathrm{SO}(q - 1,q) > (A_1)_c \times (A_1)_s > (A_1)_c \times \T > \T^2 > \T >1$$ is unrefinable and so again $\lambda(G) \leq 7$. Finally, $p=12$ and $q=4$ then the chain $$G > \mathrm{SO}(11) \times \mathrm{SO}(1,4) > \mathrm{SO}(11) \times (A_1)_c^2 > (A_1)_c^3 > (A_1)_c^2 > (A_1)_c > \T > 1$$ is unrefinable and so once more $\lambda(G) \leq 7$. 
\end{proof}


\begin{thebibliography}{30}

\bibitem{B}
 {\sc A.~Borel}, 
   {\em Linear algebraic groups}, 
Graduate Texts in Mathematics 126, 2nd Edition, Springer, New York (1969).

\bibitem{BT}
 {\sc A.~Borel and \sc J.~Tits}, 
   {\em \'El\'ements unipotents et sous-groupes paraboliques de groupes r\'eductifs. I}, 
Invent. Math. 12 (1971): 95--104.

\bibitem{BLS}
 {\sc T.C.~Burness, \sc M.W.~Liebeck and \sc A.~Shalev}, 
   {\em The length and depth of algebraic groups}, 
    to be published in Math. Zeit. (2018): 1--20.

\bibitem{BLS1}
 {\sc T.C.~Burness, \sc M.W.~Liebeck and \sc A.~Shalev}, 
   {\em The length and depth of compact Lie groups}, 
    arXiv:1805.09893 preprint (2018).

\bibitem{BLS3}
 {\sc T.C.~Burness, \sc M.W.~Liebeck and \sc A.~Shalev}, 
   {\em On the length and depth of finite groups (with an appendix by D.R. Heath-Brown)}, 
    arXiv:1802.02194 preprint (2018).

\bibitem{C1}
 {\sc R.W.~Carter}, 
   {\em Simple groups of Lie type}, 
    Wiley, New York (1972).

\bibitem{GLS}
 {\sc D.~Gorenstein, \sc R.~Lyons and \sc R.~Solomon}, 
   {\em The classification of the finite simple groups, number 3}, 
   Mathematical Surveys and Monographs 40, AMS, Providence (1997).

\bibitem{H}
 {\sc J.E.~Humphreys}, 
   {\em Linear algebraic groups}, 
    Graduate Texts in Mathematics 21, Springer, New York (1975).

\bibitem{Kn}
 {\sc A.W.~Knapp}, 
   {\em Lie Groups Beyond an Introduction}, 
Progress in Math. 140, 2nd Edition, Birkh\"{a}user, Basel (2002).

\bibitem{K}
 {\sc B.P.~Komrakov}, 
   {\em Primitive actions and the Sophus Lie problem}, 
 in: \textit{The Sophus Lie Memorial Conference} O.A. Laudal and B. Jahren eds., Scand. Univ. Press (1994): 187--269.

\bibitem{Lu}
 {\sc F.~Lubeck}, 
   {\em Small degree representations of finite Chevalley groups in defining characteristic}, 
LMS J. Comput. Math. 4 (2001): 135--169.

\bibitem{MT}
 {\sc G.~Malle and \sc D.~Testerman}, 
   {\em Linear algebraic groups and finite groups of Lie type}, 
    Cambridge Studies in Adv. Math. 133, Cambridge University Press (2011).

\bibitem{OV}
 {\sc A.L.~Onishchik and \sc E.B.~Vinberg}, 
   {\em Lie groups and algebraic groups}, 
   Springer Series in Soviet Math., New York (1988).

\bibitem{PR}
 {\sc V.~Platonov and \sc A.~Rapinchuk}, 
   {\em Algebraic groups and number theory}, 
    Pure and Applied Math. 139, Academic Press, Boston (1993).

\bibitem{S}
 {\sc R.~Steinberg}, 
   {\em Lectures on Chevalley groups}, 
    University Lecture Series 66, AMS, Providence (2016).

\bibitem{Ta}
 {\sc M.S.~Taufik}, 
   {\em On maximal subalgebras in classical real Lie algebras}, 
    Selecta Math. Sov. 6 (1987): 163--176.

\bibitem{T}
 {\sc J.~Tits}, 
   {\em Classification of algebraic semisimple groups}, 
    Proc. Sympos. Pure Math. 9 (1966): 33--62.
\end{thebibliography}
\end{document}